\documentclass[12pt,reqno]{amsart}
\usepackage{color,mathrsfs,epsfig}
\usepackage{psfrag,graphicx}
\usepackage{amssymb}

\newcommand{\ue}{u^{\varepsilon}}

\newcommand{\ee}{\varepsilon}

\newcommand{\R}{\mathbb{R}}

\newcommand{\loc}{{\rm loc}}

\newtheorem{theorem}{Theorem}[section]
\newtheorem{lemma}[theorem]{Lemma}

\theoremstyle{definition}
\newtheorem{definition}[theorem]{Definition}

\theoremstyle{remark}
\newtheorem{remark}[theorem]{Remark}
\newtheorem{example}[theorem]{Example}

\newtheorem{hyp}{Hypothesis}

\numberwithin{equation}{section}

\begin{document}
\bibliographystyle{plain}
\title[Viscous approximation via ODE analysis]{Notes on the study of the viscous approximation of hyperbolic problems via ODE analysis}
\author{Laura V.~Spinolo}
\address{Centro De Giorgi, Collegio Puteano, Scuola Normale Superiore,
Piazza dei Ca\-va\-lie\-ri 3, 56126 Pisa, Italy}
\email{laura.spinolo@sns.it}
\begin{abstract} These notes describe some applications of the analysis of ordinary differential equations to the study of the viscous approximation of conservation laws in one space dimension. The exposition mostly focuses on the analysis of invariant manifolds like the center manifold and the stable manifold. The last section addresses a more specific issue and describes a possible way of extending the notions of center and stable manifold to some classes of singular ordinary differential equations arising in the study of the Navier-Stokes equation in one space variable. 
\end{abstract}

\maketitle

\section{ Introduction}
\label{s:out}
In recent years, techniques coming from the study of ordinary differential equations (ODEs) have been applied to the analysis of systems of conservation laws in one space dimension, 
\begin{equation}
\label{e:cl}
    u_t + f(u)_x =0 \quad \textrm{where $u (t, x) \in \R^N $ and $(t, x) \in [0, + \infty [ \times \R$}.
\end{equation}
In the previous expression, ${}_t$ and ${}_x$ denote the partial derivatives with respect to the variables $t$ and $x$ respectively and the flux $f: \R^N \to \R^N$ is a function of class $\mathcal C^2$. In view of physical considerations  (see e.g. Dafermos~\cite{Daf:book}) it is of great interest considering the second order approximation 
\begin{equation}
\label{e:va}
    \ue_t + f(\ue)_x = \ee \big( B( \ue ) \ue_x \big)_x, 
\end{equation}
where $B$ is an $N \times N$ ``viscosity matrix" satisfying suitable conditions and $\ee$ is a parameter, $\ee \to 0^+$. In particular, the Navier-Stokes equation in one space variable can be written in the form~\eqref{e:va} and by setting $\ee =0$ it formally reduces to the Euler equation. For the time being, we focus on the Cauchy problem obtained by coupling~\eqref{e:va} with the initial datum
\begin{equation}
\label{e:cauchy}
    u(0, x) = u_0 (x). 
\end{equation}
Establishing general results on the convergence of the family of functions $\ue$ satisfying~\eqref{e:va},~\eqref{e:cauchy} is still a mayor open problem, but partial results have been achieved and ODE techniques have often played an important role: see again Dafermos~\cite{Daf:book} for a complete overview and a list of references. Let us just mention the work by Bianchini and Bressan~\cite{BiaBrevv}, which established convergence (under suitable assumptions) of the family of functions $\ue$ solving~\eqref{e:va},~\eqref{e:cauchy} in the case of the ``artificial viscosity" $B(\ue) \equiv I$ (identity matrix). The proof in~\cite{BiaBrevv} employs a special decomposition of the gradient $\ue_x$ obtained by relying on ODE analysis. 

These notes aim at providing an informal overview on some of the ODE techniques that are more frequently applied to the study of the limit $\ee \to 0^+$ of~\eqref{e:va}. The exposition is addressed to non-experts, so notions are usually first introduced in a simplified context and then discussed in more general situations. Also, in most cases these notes only provide an heuristic idea of the proof of the the results that are introduced, and refer to books or to original research papers for a more rigorous discussion. The discussion has, of course, no sake of completeness: in particular, to simplify the exposition very few references are provided. For a more satisfactory bibliography, one can refer to the books by Dafermos~\cite{Daf:book} and by Serre~\cite{Serre:book} (conservation laws and their viscous approximation). An extremely rich exposition of the analysis of invariant manifolds for ODEs is in the book by Katok and Hasselblatt~\cite{KatokHass} (see also the book by Perko~\cite{Perko}). 

To simplify the exposition, most of the following sections (all but the last one) focus on the case of the ``artificial viscosity" $B(\ue) \equiv I$, so that~\eqref{e:va} reduces to
\begin{equation}
\label{e:av}
     \ue_t + f(\ue)_x = \ee \ue_{xx}.
\end{equation}
 However, many of the considerations extend to more general cases (sometimes the extension is not straightforward, though). The details concerning this extension can be found in the original research papers that are quoted in the following sections. 

Also, most of the analysis described in the following sections apply to the study of the non-conservative case 
\begin{equation}
\label{e:nc}
    \ue_t + A(\ue, \ee \ue_x) \ue_x = \ee B(\ue) \ue_{xx},
\end{equation}
where $A$ is a suitable $N \times N$ matrices. Note that~\eqref{e:va} can be written in the form~\eqref{e:nc} by setting 
$
    A (\ue, \ee \ue_x) = D f (\ue) - \ee B(\ue)_x, 
$
where $Df$ denotes the Jacobian matrix of the flux $f$. Again, the details concerning this extension are available in the original research papers quoted in the following sections.

The exposition is organized as follows: Section~\ref{s:mot} informally describes the connection between the analysis of the limit $\ee \to 0^+$ of the viscous approximation~\eqref{e:va} and the study of ODEs by introducing the notions of \emph{traveling waves} and \emph{boundary layers}. Section~\ref{s:cm} deals with the Center Manifold Theorem and Section~\ref{s:cm:app} mentions some applications to the study of the viscous approximation. Section~\ref{s:stable} deals with the Stable Manifold Theorem, while the last section of these notes (Section~\ref{s:sode}) is more technical and describes a possible way of extending the definition of center and stable manifold to a class of singular ordinary differential equations arising in the study of the viscous approximation defined by the Navier-Stokes equation in one space variable. 

\subsubsection*{Acknowledgments}
These notes were originally prepared for a course given at the {\em Seventh Meeting on
Hyperbolic Conservation Laws and Fluid Dynamics:
Recent Results and Research Perspectives}. The conference was held at SISSA, Trieste and was organized by Fabio Ancona, Stefano Bianchini, Rinaldo M. Colombo and Andrea Marson. The author warmly thank them for the kind invitation. 
\section{ Some motivations}
\label{s:mot}
This section aims at discussing some links between the study of the viscous approximation of a system of conservation laws in one space dimension and the ODE analysis by informally illustrating the notions of traveling waves and boundary layers.

\subsection{Traveling waves}
\label{sus:tw}
Let us start with some heuristic considerations. For every $\ee >0$, \eqref{e:av} is a parabolic equation and one expects a strong regularizing effect. Conversely, it known that, even if the initial datum $u_0$ is smooth, nevertheless in general a classical solution (i.e., a solution of class $\mathcal C^1$)  of the Cauchy problem~\eqref{e:cl},~\eqref{e:cauchy} may break down in finite time: remarkably, this may happen even in the scalar case $N=1$ (see e.g. Dafermos~\cite[Chapter 4.2]{Daf:book}). Hence, when studying the limit $\ee \to 0^+$ one handles a family of very regular functions which is expected to converge to something which has a much weaker regularity. Traveling waves are often regarded as useful tools to handle this behavior and to describe the formation of singularities in the limit.  
 
 To fix the notations, let us introduce the following definition.
 \label{sss:tw}
  \begin{definition}
 \label{d:tw}
           Let $u^+, u^- \in \R^N$ and $\sigma \in \R$ be given. 
          The function $U: \R \to \R^N$ is a traveling wave for~\eqref{e:va} 
           joining the states $u^+, u^- \in \R^N$ and having speed $\sigma$ if 
           $U(y)$ satisfies 
            \begin{equation}
            \label{e:dtw}
                       U'' = \big[ f(U) - \sigma U \big]' 
             \end{equation}
             and 
             \begin{equation}
             \label{e:ltw}          
                       \lim_{ y \to - \infty} = u^- \qquad \lim_{ y \to + \infty} = u^+ 
            \end{equation} 
 \end{definition}
In the previous expression, $U'$ and $U''$ denote respectively the first and the second derivative of $U$ with respect to the variable $y$. Some observations are here in order. First, the link between~\eqref{e:dtw} and~\eqref{e:av} is the following: if $U$ solves~\eqref{e:dtw}, then one can verify by direct check that a solution of~\eqref{e:av} is obtained by setting 
\begin{equation}
\label{e:dtw2}
     \ue (t, x) : = U \left( \frac{x - \sigma t}{\ee} \right).
\end{equation}
Also, the function $U$ is a solution of the ordinary differential equation~\eqref{e:dtw}, hence it is of class $\mathcal C^2$ and the same regularity is inherited by $\ue$. However, by taking the pointwise limit $\ee \to 0^+$ of $\ue$ we get 
$$
    u(t, x) : = \lim_{\ee \to 0^+} \ue (t, x) =  \lim_{\ee \to 0^+} 
    U \left( \frac{x - \sigma t}{\ee} \right) =
    \left\{
    \begin{array}{lll}
                \lim_{y \to - \infty} U(y) \qquad \text{if $ x - \sigma t < 0$} \\
                \\
                \lim_{y \to + \infty} U(y) \qquad \text{if $ x - \sigma t > 0$} \\           
    \end{array}
    \right.
$$
One can then verify that the pointwise limit $u(t, x)$ provides a distributional solution of~\eqref{e:cl}.

Summing up, we have that, if conditions~\eqref{e:ltw} hold, then the solution of the parabolic equation defined by~\eqref{e:dtw2} is regular for every $\ee >0$, but converges pointwise to a discontinuous distributional solution of the conservation law. This is why traveling waves are often used to study the non trivial behavior mentioned above, namely the appearance of singularities when passing from the vanishing viscosity approximation to the hyperbolic limit. 

Let us now why the ODE analysis comes into play: the following considerations are quite informal, for a more precise and rigorous exposition see for example Bianchini and Bressan~\cite[Sections 3, 4]{BiaBrevv}. Let us consider the equation~\eqref{e:dtw}: usually, the values $u^+$, $u^-$ and $\sigma$ are not prescribed. Conversely, one wants to understand for which values of $u^+$, $u^-$ and $\sigma$ there is a solution of~\eqref{e:dtw} satisfying 
           $$
                 \lim_{ y \to - \infty} = u^- \qquad \lim_{ y \to + \infty} = u^+ . 
           $$                     
     To select a unique solution of~\eqref{e:dtw}, one needs to assign $2N +1$ parameters: for example, one can assign a Cauchy datum for $U$ and for $U'$ and the value of $\sigma$ (remember that $U$ takes values in $\R^N$). As a matter of fact, in many situations one does not manage to assign so many conditions and apparently one ends up with an underdetermined problem. By assigning $2N +1$ parameters, however, one is neglecting an important information: the values of $u^+$ and $u^-$ are not prescribed, but nevertheless both the limits $\lim_{x \to \pm} U(y)$ exist and are finite. In particular, this implies that $U$ is bounded on the whole $\R$. In many situations, one actually has even more information: given a value $\bar u \in \R^N$, one looks for values of $u^+$, $u^-$ and $\sigma$       
satisfying the following property. 
\begin{equation}
\label{e:cmc}
    \textrm{Both $u^+$ and $u^-$ are close to $\bar u$, $\sigma$ is close to $\lambda_i (\bar u)$}.
\end{equation}
Here $\lambda_i (\bar u)$ denotes  the $i$-th eigenvalue of the Jacobian matrix $Df(\bar u)$.  
As we see in Section~\ref{s:cm}, the Center Manifold Theorem allows to exploit the information~\eqref{e:cmc} to reduce the number of conditions one needs to impose on~\eqref{e:dtw} to select a unique solution. 
\subsection{Boundary layers}
\label{sus:bl}
Boundary layer phenomena are observed when studying the viscous approximation of initial boundary value problems.

To illustrate the heart of the matter, let us focus on a specific situation. Consider the family
of problems 
$$
\left\{
\begin{array}{lll}
            \ue_t + f(\ue)_x = \ee \ue_{xx} \\
            \ue (t, 0) = u_b (t) \qquad \ue (0, x) = u_0 (x). 
\end{array}
\right.
$$
As established in Ancona and Bianchini~\cite{AnBiapro}, under suitable assumptions on the data $u_b$ and $u_0$, for every $t \ge 0$ the family $\ue (t, \cdot)$ converges in $L^1_{\loc}$ to $u(t, \cdot)$, a distributional solution of~\eqref{e:cl} satisfying the following condition:
for every $t\ge 0$, $\mathrm{TotVar}_x u(t, \cdot) < + \infty$. The point where boundary layers come into play is the fact that, in general, the boundary condition is lost, namely
(note that the limit below exists since the total variation is bounded) 
\begin{equation}
\label{e:trace}
    \lim_{x \to 0^+} u(t, x) \neq u_b (t). 
\end{equation}
See also the previous work by Gisclon~\cite{Gisclon:etudes} for the proof of the local in time convergence for more general viscous approximations~\eqref{e:va}. 

To study the ``loss" of boundary condition~\eqref{e:trace}, the key point is the analysis of the steady solutions of~\eqref{e:av}, satisfying $f(u)_x = \ee u_{xx}$. To understand why steady solutions are able to capture this behavior (the loss of boundary condition), let us start by considering a toy model. 
\begin{example}
\label{ex:bl}
Consider the vanishing viscosity approximation of a scalar, linear conservation law
\begin{equation}
\label{e:ex}
    \ue_t + a \ue_x = \ee \ue_{xx} \quad a < 0
\end{equation}
and let us focus on the steady solutions $\ue (t, x) = U_{\ee}(x)$, obtained by solving the linear ODE
\begin{equation}
\label{e:steady}
    a U_{\ee}' = \ee U_{\ee}''.
\end{equation}
Also, let us impose on~\eqref{e:steady} the conditions 
$$
    U_{\ee}(0) = u_b  \qquad \lim_{x \to + \infty} U_{\ee}(x) = u_0. 
$$
Then by performing explicit computations we get 
$
    U_{\ee}(x ) = \big[ u_b - u_0 \big] e^{a x/ \ee } + u_0. 
$
Hence, by letting $\ee \to 0^+$ and, by exploiting $a <0$, we obtain that, for every $x>0$, $U_{\ee}(x)$ converges to $U(x) = u_0$ and hence the boundary condition $U_{\ee}(x) = u_b$ is ``lost" since $\lim_{x \to 0^+} U(x) = u_0$. 
\end{example}
Some observations are here in order: first, by introducing the change of variables  $(t, x) \mapsto (\ee t, \ee x)$ one can transform the equation $f(\ue)_x = \ee \ue_{xx}$ into $f(u)_x = u_{xx}$. Also, assume that we want to study the ``loss of boundary condition" for a linear system of conservation laws 
$$
    u_t + A u_x =  u_{xx} \quad \text{where $u \in \R^N$}.
$$
Here, $A$ is an $N \times N$ matrix and we assume that it admits $N$ real and distinct eigenvalues.
We want to proceed as in Example~\ref{ex:bl} by focusing on the steady solutions $u(x, t) = U(x)$ and assigning the conditions $U(0) = u_b$ and $\lim_{x \to + \infty} U(x) = u_0$. However, it turns out that these conditions are compatible if and only if the vector $u_b - u_0$ belongs to the subspace of $\R^N$ generated by the eigenvectors of $A$ associated to \emph{strictly negative} eigenvalues.  Let us now consider the general case: the second order approximation of a non linear system of conservation laws 
$$
     u_t + f ( u)_x = u_{xx} \quad \text{where $u \in \R^N$}.
$$
To extend the argument in Example~\ref{ex:bl} , we have to study the system 
\begin{equation}
\label{e:d:bl}
          \left\{
    \begin{array}{lll}
                \big[ f(U) \big]' = U '' \\  \\
                U( 0) = u_b \qquad \lim_{x \to + \infty} U(x) = u_0. 
    \end{array}
    \right.
\end{equation}
In view of the previous considerations it is thus natural to ask ourself the question: what are the values of $u_b$ and $u_0$ ensuring that~\eqref{e:d:bl} admits a solution?
Note that such a query has some similarities with the problem we discussed at the end of Section~\ref{sss:tw}, namely we are concerned with an ODE, we want to prescribe some asymptotic behavior and we ask ourselves what are the data compatible with such an asymptotic behavior. 

As we will see, an answer to the above question can be obtained by relying on the Stable Manifold Theorem, discussed in Section~\ref{s:stable}.  

Finally, to fix again the notations, we introduce the following definition. 
 \begin{definition}
\label{d:bl} Let the states $u_b$, $u_0 \in \R^N$ be given. A boundary layer for (2.1) connecting
the states $u_b$ and $u_0$ is a solution of system~\eqref{e:d:bl}.
\end{definition}


\section{ The Center manifold Theorem}
\label{s:cm}
In this section we are concerned with studying the asymptotic behavior of the solutions of  the ordinary differential equation
$$
    V' = G(V) \quad \text{where $V \in \R^d$}
$$
in a neighborhood of an equilibrium point $\bar V$. Without any loss of generality we can assume $\bar V = \vec 0$, namely $G(\vec 0) = \vec 0$. In Section~\ref{s:cm:app} we will then apply our analysis to the study of the system~\eqref{e:dtw}. 

Let us start by considering two simple examples where $G$ is linear. 
\begin{example}
\label{ex:cm:lin}
Assume that $V=(x, y, w, z)$ and that 
\begin{equation}
\label{e:ex:lin}
         \left(
         \begin{array}{cccc}
                   x' \\ y' \\ w' \\ z ' \\ 
         \end{array}
         \right) = 
          \left(
         \begin{array}{cccc}
                   2 & 0 & 0 & 0 \\ 0 & -1  & 0 & 0  \\ 0 & 0  & 0 & - 3 \\   0 & 0  & 3 & 0 \\ 
         \end{array}
         \right) 
          \left(
         \begin{array}{cccc}
                   x \\ y \\ w \\ z  \\ 
         \end{array}
         \right) 
\end{equation}
By solving the ODE explicitly we obtain 
$$
     \left(
         \begin{array}{cccc}
                   x (t)  \\  y (t)  \\ w (t) \\ z (t)  \\ 
         \end{array}
         \right) = 
     \left(
         \begin{array}{cccc}
                   \displaystyle{e^{2 t} x(0)}   \\  
                   \displaystyle{ e^{-t} y(0)}    \\  
                   w(0) \cos 3 t  - z (0) \sin 3 t         \\ 
                   z (0) \cos 3 t   + w(0) \sin 3t   \\ 
         \end{array}
         \right)     
$$ 
Note that the eigenvalues of the above matrix are $2$, $-1$ and $\pm 3 i$. If the initial datum belongs to the eigenspace associated to the eigenvalue with \emph{positive} real part, namely if the initial datum is in the form
$$
    \Big( x(0), 0, 0, 0\Big),
$$ 
then the solution blows up exponentially fast when $x \to + \infty$, while when $x \to - \infty$ it converges with exponential speed to the equilibrium point. Conversely, assume that the initial datum belongs to the eigenspace associated to the eigenvalue with \emph{negative} real part, namely the initial datum is in the form
$$
    \Big( 0, y(0), 0, 0 \Big).
$$ 
Then the solution converges exponentially fast to the equilibrium for $x \to + \infty$, while it blows up exponentially fast when $x \to - \infty$. Finally, assume that the initial datum belongs to the eigenspace associated with the eigenvalues with \emph{zero} real part, namely it is in the form
$$
    \Big( 0, 0,  w(0), z(0) \Big).
$$   
Then the solution does not converge exponentially neither at $+ \infty$ nor at $ -\infty$.
However, the solution is \emph{globally bounded} on the whole real line. 
\end{example}
The second example focuses on a linear system where, however, the matrix is not diagonalizable.         
\begin{example}
\label{ex:cm2}
          Assume that $V = (x, y)$ and that 
          \begin{equation}
          \label{e:ex:cm2}
                    \left(
         \begin{array}{cc}
                   x' \\ y' \\  
         \end{array}
         \right) = 
          \left(
         \begin{array}{cc}
                   0 & 1  \\ 0 & 0  \\ 
          \end{array}
         \right) 
          \left(
         \begin{array}{cc}
                   x \\ y \\ 
         \end{array}
         \right) 
          \end{equation}
      The above matrix has eigenvalue $0$ with multiplicity $2$ and hence the whole $\R^2$ is the eigenspace associated with eigenvalues with zero real part. By solving the equation explicitly we get
     $$
         \left(
         \begin{array}{cc}
                   x \\ y \\ 
         \end{array}
         \right) 
         =
            \left(
         \begin{array}{cc}
                   x(0) + y(0) t \\ y (0) \\ 
         \end{array}
         \right) .
     $$ 
     In this case the solution is not bounded on the whole $\R$ unless $y(0) =0$. However, it does not blow up exponentially neither at $x \to + \infty$ nor at $ x \to - \infty$.  
\end{example}

\subsection{The center manifold theorem}
\label{sus:cm}
\subsubsection{Statement of the theorem}
\label{sss:state}
Let us first introduce some notations: $B (\vec 0, \delta)$ is the ball with radius $\delta$ and center at $\vec 0$ in $\R^d$. Also, $DG(\vec 0)$ is the Jacobian matrix of $G$
         computed at 
         $V = \vec 0$ and $V^c$ will be its center subspace, namely
 \begin{equation}
         \label{e:vc}
             V^c = Z_1 \oplus Z_2 \oplus \dots Z_{n_c},
   \end{equation}        
   where each of the $Z_i$ is an eigenspace associated to an eigenvalue $\lambda_i$ with zero real part, 
   $$
             Z_i = \big\{ \vec v \; \text{such that $[DG(\vec 0) - \lambda_i I]^ k \vec v = \vec 0$ for some $k \leq d $} \big\} 
             \quad \mathrm{Re} \lambda_i =0.         
    $$
We are now ready to introduce the main result of this section: 
\begin{theorem}
\label{t:cm}
          Consider the first order ODE
          \begin{equation}
          \label{e:cm:ode}
              V' = G(V), \quad V \in \R^d,
         \end{equation}
         where $G: \R^d \to \R^d$ is a $\mathcal C^2$ function satisfying 
         $G (\vec 0) = \vec 0$. Assume that the center space $V^c$ defined by~\eqref{e:vc} is 
         non trivial, namely $ V^c  \neq \big\{ \vec 0 \big\}$.
                 Then there exists a constant $\delta >0$ small enough and a continuously differentiable  
         center manifold $\mathcal M^c$
         satisfying the following conditions:
         \begin{enumerate}
         \item $\mathcal M^c$ is parameterized by a function
         $$
             \phi_c : B(\vec 0, \delta) \cap V^c \to \R^d.
         $$
          Also, $\mathcal M^c$ is tangent to $V^c$ at the origin.
         \item $\mathcal M^c$ is locally invariant for the ODE~\eqref{e:cm:ode}. Namely,  
         if $V_0 \in \mathcal M^c$, then the solution of the Cauchy problem
         $$
         \left\{
         \begin{array}{lll}
                   V' = G(V) \\
                   V(0) = V_0     
         \end{array}
         \right.
         $$
         belongs to $\mathcal M^c$ if $|t|$ is small enough. 
         \item let $V(t)$ be a solution of~\eqref{e:cm:ode} 
         such that $V(t) \in B(\delta, \vec 0)$ for every $t \in \R$. 
         Then 
         $$
              V(t) \in \mathcal M^c \quad \text{ for every $t \in \R$}. 
         $$     
        \end{enumerate}                  
\end{theorem}
An important remark is that the center manifold~\eqref{e:cm:ode} about the equilibrium point $\vec 0$ is not unique (see for example the lecture notes by Bressan~\cite{Bre:note} for an explicit example). In other words, suitable examples show that in general there is more than one manifold satisfying conditions (1)..(4) in the statement of the theorem. We will come back to this point in the proof of the theorem, when it will be clear where this lack of uniqueness come from.

The proof of Theorem~\ref{t:cm} discussed in here is the same as in the notes by Bressan~\cite{Bre:note}. For a more general viewpoint, one can refer to the book by Katok and Hasselblatt~\cite{KatokHass}. 
\subsection{Proof of Theorem~\ref{t:cm}}
\label{sus:proof}
The proof is divided in several steps, and only some of them are provided below (see~\cite{Bre:note} for the complete argument). \\
{\sc Step 1: heuristic.} The heuristic idea is that the manifold $\mathcal M^c$ should play the role that in the linear case is played by the center space. In exploiting this idea one should keep in mind that, as Example~\ref{ex:cm2} shows, the solutions of a linear system lying 
on the center space in general are not bounded on the whole real line: the only requirement we can reasonably impose is that they did not blow up exponentially fast neither 
at $+ \infty$ nor at $- \infty$. 

To make an extension to the general non linear case we have to handle some technical difficulties. First of all, we need to introduce a localization argument. Also, because of the non linearity we cannot require that the solutions lying on the center manifold do not blow exponentially neither 
at $+ \infty$ nor at $- \infty$, but only that they are controlled by $\exp (\eta |t|)$, where $\eta$ is a small positive constant. \\
{\sc Introduction of a cut-off function} Let 
$$
     p(V) = G(V) - DG(0) V
$$
be the second order term in the expansion of $G$. We then have 
\begin{equation}
\label{e:secondo}
          |p(V) | \leq C | V|^2 \quad \textrm{and} \quad |D p (V)| \leq C V
\end{equation}
for a suitable constant $C$ and for $V$ in a small enough neighborhood of $\vec 0$.

We want to reduce to the case $p$ has compact support and small $\mathcal C^1$ norm. We thus introduce a smooth, even cut-off function $\rho: \R \to \R$ satisfying 
$$
     0 \leq \rho(t) \leq 1 \; \text{for every $t \in \R$} \qquad \text{ and} \qquad   
    \rho (t) =
    \left\{
    \begin{array}{lll}
               1 & \text{if $|t| \leq  1$} \\
               0 & \text{if $|t| \ge 2$}\\
    \end{array}
    \right.
$$
We then define 
$$
    p_{\delta} (V) : = p(V) \, \rho ( |V| / \delta),
$$
thus obtaining that the support of $p_{\delta}$ is contained in $B (\vec 0, 2 \delta)$. In the previous expression, $\delta$ denotes a small constant whose exact value is determined in the following. By relying on~\eqref{e:secondo} we also get  
$$
     \|p_{\delta} \|_{\mathcal C^0} 
     \leq \sup_{ X \in B (\vec 0, 2 \delta)} |p (V) | 
     \leq 4 C \delta^2 
$$
and 
$$
    | D p_\delta (V) |= |D p (V)| \rho + |p(V)| \,  |\rho '|   \, \frac{1}{\delta} 
    \leq \tilde C \delta. 
$$
By combing the previous observations we can make the $\mathcal C^1$ norm arbitrarily small by choosing $\delta$ small enough.

In the following, we will study the solutions of the ODE
\begin{equation}
\label{e:cutoff}
    V' =  DG(0) V + p_\delta (V).
\end{equation}
By construction, 
$$
    DG(0) V + p_\delta (V) = G(V) \quad \text{when $|V| \leq \delta$}
$$
and hence any claim holding for the solutions of~\eqref{e:cutoff} holds for the solutions of the original equation 
$
    { V' = G(V)}
$ confined in $B(\vec 0, \delta)$. \\
{\sc Step 3: definition of a fixed point problem}
We first have to introduce some further notations: for simplicity, in the following we denote by $A$ the Jacobian matrix $DG (\vec 0)$. Also, the stable and unstable space of $A$ are the subspaces of $\R^d$ associated to the eigenvalues with negative and positive real part respectively. Namely, 
\begin{equation}
\label{e:stableun}
     V^s = N_1 \oplus N_2 \oplus \dots N_{n_s} \qquad 
     V^u = P_1 \oplus P_2 \oplus \dots P_{n_u} 
  \end{equation}
where the $N_i$ and the $P_i$ are eigenspaces associated with eigenvalues with negative and positive real part respectively:
$$
   N_i= \big\{ \vec v \in \R^d: \;  \big( A - \lambda_i I \big)^k v = \vec 0 \; 
    \textrm{for some $k \leq d$} \} \quad \mathrm{Re} \lambda_i <  0 $$
and 
$$
   P_i= \big\{ \vec v \in \R^d: \;  \big( A - \lambda_i I \big)^k v = \vec 0 \; 
    \textrm{for some $k \leq d$} \} \quad \mathrm{Re} \lambda_i >  0. $$In particular, since 
         $$
             \R^d = V^s \oplus V^u \oplus V^c,
         $$
         then any $x \in \R^d$ satisfies $x = \pi_s x + \pi_u x+ \pi_c x$, where    
      $$
             \pi_s: \R^d \to V^s \qquad \pi_u: \R^d \to V^u \qquad \pi_c: \R^d \to V^c
         $$
         are the projections onto $V^s, \; V^u$ and $V^c$ respectively.  
         
Note that the ODE
         $$
             V' = A V + p_\delta (V)
         $$
         can be then written as the coupling of the ODEs
         $$
            \left\{
            \begin{array}{lllll}
                       \displaystyle{ \big( \pi_s V \big)' =   A \pi_s V + \pi_s p_\delta (V) }\\
                       \\
                       \displaystyle{ \big( \pi_u V \big)' =  A \pi_u  V + \pi_u p_\delta (V) } \\
                       \\
                       \displaystyle{ \big( \pi_c V \big)' =   A \pi_c V + \pi_c p_\delta (V) }. \\ 
            \end{array}
            \right.
         $$
         We have exploited the following equalities:
         $$
             \pi_s A V = A \pi_s V   \quad   \pi_u A V = A \pi_u  V  \quad 
               \pi_c A V = A \pi_c V \quad \textrm{for every $V \in \R^d$}
         $$ 
         We solve each equation separately and by applying the 
         ``variation of constants" formula we obtain that the components $\pi_s V(t)$,
           $\pi_u V(t)$ and  $\pi_c V(t)$ satisfy 
         \begin{equation}
         \label{e:ode:comp}
             \left\{
             \begin{array}{lllll}
            \displaystyle{ \pi_s V (t) = e^{  A (t-t_s)}  
            \pi_s V (t_s) + \int_{t_s}^t  e^{ A (t-s) }
             \pi_s p_\delta (V (s)) ds}  \\
             \\
            \displaystyle{\pi_u V (t) = e^{A (t-t_u)}  
            \pi_u V (t_u) + \int_{t_u}^t  e^{ A (t-s) }
             \pi_u p_\delta (V (s)) ds}  \\
             \\ 
            \displaystyle{\pi_c V (t) = e^{A (t-t_c)}  
            \pi_c V (t_c) + \int_{t_c}^t  e^{A (t-s) }
             \pi_c p_\delta (V (s)) ds}.
            \end{array}
            \right.
         \end{equation}
         In the previous expression, $t_s$, $t_u$ and $t_c$ 
         denote real values which will be assigned in the following. Also, we have exploited the standard notation 
         $$
               \exp ( A t )  = \sum_{n =0}^{\infty} \frac{t^n A^n}{n!}.
         $$
         Let $\lambda_i$ denote as before an eigenvalue of $A$, 
         then we set
         $$
             \beta_+ : = \min \{\mathrm{Re} \lambda_i: \; \mathrm{Re} \lambda_i >0 \}
         $$
         and 
         $$
            \beta_- : = \min \{ - \mathrm{Re} \lambda_i: \; \mathrm{Re} \lambda_i < 0 \}.
         $$
         There exists some constant $C>0$ such that, for every $\vec v \in \R^d$, we have 
         \begin{equation}
         \label{e:exp}
             |e^{A t} \pi_u \vec v | \leq C e^{\beta_+ t / 2 } |\vec v | \; \text{for every $t < 0$} \; \text{and} \; 
              |e^{A t} \pi_s \vec v | \leq C e^{- \beta_- t/ 2 } |\vec v | \; \text{for every $ t > 0$}. 
         \end{equation}
          We now fix a constant positive constant $\eta >0$ satisfying 
$$
    \eta <\min\left\{  \frac{ \beta_+ }{2} , \;   \frac{\beta_- }{2} \right\}
$$    
and set
$$
    Y_{\eta} : = \big\{ V \in \mathcal C^0 ( ]- \infty, + \infty [ ; \R^d \big\} 
    \; \text{such that} \; \sup_{t \in \R} |V (t) | e^{- \eta |t|} < + \infty   \big\}. 
$$
Then $Y_{\eta}$ equipped with the norm
$$
     \| V \|_\eta : = \sup_{t \in \R} |V (t) | e^{- \eta | t |} 
$$
is a Banach space. Note that 
$$
    |V(t) | \leq e^{ \eta |t|}   \| V \|_\eta  \quad  \text{for every $t \in \R$, $V \in Y_\eta$}.
$$
Loosely speaking, the space $Y_\eta$ contains all the function which ``do not blow up with too fast exponential speed" neither at $+\infty$ nor at $- \infty$. 

The goal is now defining the map $\phi_c$ by solving a suitable fixed point problem in $Y_{\eta}$. First, we observe that, if the function $V \in Y_{\eta}$, then by relying on~\eqref{e:exp}  we get that, for any $t_s \leq t$,
$$
    |e^{A (t-t_s)}  
            \pi_s V (t_s) | \leq e^{- \beta_- (t - t_s) / 2} e^{- \eta  t_s } \| V \|_\eta =
    e^{ - \beta_- t / 2 } \| V \|_\eta e^{( -  \eta + \beta_- / 2 )  t_s }        
$$
and hence
$$
     \lim_{t_s  \to - \infty} |e^{ A (t-t_s)}  
            \pi_s V (t_s) | = 0 \quad \text{for any given $t \in \R$}. 
$$
In the same way one can show that, if the function $V \in Y_\eta$, then 
$$
     \lim_{t_u  \to + \infty} |e^{ A (t-t_u)}  
            \pi_u V (t_u) | = 0 \quad \text{for any given $t \in \R$}. 
$$
We now go back to formula~\eqref{e:ode:comp}, we set $t_c =0$, $x_c = \pi_c V (0)$ and we let $t_s \to - \infty$ and $t_u \to + \infty$, eventually obtaining 
\begin{equation}
\label{e:fp}
    \begin{split}
    V(t) = 
    &      \pi_c V(t) + \pi_s V(t) + \pi_u V(t)  =   e^{A t } x_c 
             + \int_{0}^t  e^{A (t-s) }
             \pi_c p_\delta (V (s)) ds \\ &+  \int_{- \infty}^t  e^{ A (t-s) } 
             \pi_s p_\delta (V (s)) ds 
            +  \int_{+ \infty}^t  e^{A (t-s) }
             \pi_u p_\delta (V (s)) ds  
          \\
    \end{split}         
\end{equation}
for every $V \in Y_\eta$. \\
{\sc Further steps in the proof of the Center Manifold Theorem}
\begin{itemize}
\item
By applying the Contraction Map Theorem, one can show that, for any $x_c  \in V^c$, the fixed point problem~\eqref{e:fp} admits a unique solution
$V \in Y_\eta$. To verify that the hypotheses of the Contraction Map Theorem are satisfied one has to exploit that $p_\delta$ is a function with compact support and small $\mathcal C^1$ norm.
\item The map
$$
    \phi_c : V^c \to \R^d 
$$
which parameterizes the center manifold is then defined by setting 
\begin{equation}
\label{e:phic}
         \phi_c (x_c ) : = V(0).
\end{equation}         
\item One can then prove that such a map is continuously differentiable and that the tangent space at the origin is $V^c$. 
\item By construction, the map $\phi_c$ satisfies the following properties: first, for any $ x_c \in V^c$, the solution of the Cauchy problem
\begin{equation}
\label{e:fake}
\left\{
\begin{array}{lll}
            V' = AV + p_\delta (V) \\
            \\
            V(0) = \phi_c (x_c)
\end{array}
\right.
\end{equation}
belongs to $Y_\eta$, and hence it ``does not blow up with fast exponential speed" neither at $+ \infty$ nor at $-\infty$.  
Conversely, any solution of~\eqref{e:fake} verifying 
$$
   \sup_t |V(t)| e^{- \eta |t| } < + \infty
$$
satisfies $V(0) = \phi_c (\bar x_c)$ for some $\bar x_c \in V^c$.
\item We set $\mathcal M^c : = \phi_c ( V_c)$ and we eventually come back to the original equation
$
    V'  = G(V).
$
As pointed out before, $G(V) = DG(0) V + p_\delta (V) $ when $|V| < \delta$. Consequently, any claim concerning the solution of $V'=   DG(0) V + p_\delta (V) $ can be reformulated in a claim concerning the solutions of $V' = G(V)$ which are confined in $B(\vec 0, \delta)$ for every $t \in \R$. By relying on this observation one concludes the proof of Theorem~\ref{t:cm}. 
\end{itemize}
{\sc A remark concerning non-uniqueness}
As pointed out before, in general the center manifold of~\eqref{e:cm:ode} about a given equilibrium point is not unique. The reason is that in the proof we have some freedom in choosing the cut-off function $\rho$: hence, both the map $\phi_c$ defined as in~\eqref{e:phic} and the manifold $\mathcal M^c$ depend on the choice of $\rho$ and are not unique.

\section{ Some applications of the Center Manifold Theorem to the analysis of traveling waves}
\label{s:cm:app}
The Center Manifold Theorem has been applied by several authors to the analysis traveling waves: see Dafermos~\cite{Daf:book} for a list of references. In this section we will discuss a possible approach, due to Bianchini and Bressan~\cite{BiaBrevv} and in the second part we will deal with applications to the study of the vanishing viscosity approximation $\ue_t + f(\ue)_x = \ee \ue_{xx}$.  
\subsection{A center manifold of traveling waves}
\label{sus:cm:come}
Let us now go back to the problem we discussed at the end of Section~\ref{sus:tw}: we are given the traveling wave equation
\begin{equation}
\label{e:twdue}
         U'' = \big[ f (U) - \sigma U \big]' \quad \text{where $U \in \R^N$}
\end{equation}
and a value $\bar u \in \R^N$. We want find solutions $U$ admitting finite limits both at $x \to \pm  \infty$,
$
  \lim_{ y \to \pm \infty} = u^{\pm} 
$
and satisfying the following properties: first, for every $y \in \R$, $U(y)$ belongs to a small neighborhood of $\bar u$. Also, $\sigma$ is close to $\lambda_i (\bar u)$, where $\lambda_i (\bar u)$ is a given eigenvalue of the Jacobian matrix $D f (\bar u)$. 
More precisely, we want to exploit the information we have on $U$ to reduce the number of conditions one has to assign on~\eqref{e:twdue} to select a unique solution. 

By setting 
$
    p : = U'      
$
we get that~\eqref{e:twdue} can be written as the first order system
\begin{equation}
\label{e:tre}
         \left\{
         \begin{array}{lll}
                    U' = p \\
                    p' = \big[ Df (U) - \sigma I \big] p \\
                    \sigma' =0 .
         \end{array}
         \right.
\end{equation}
We now apply the Center Manifold Theorem: first, we observe that the point 
$
    {\bar V= \big( \bar U, \vec 0,  \lambda_i (\bar U) \big)}
$
is an equilibrium for~\eqref{e:tre}. Hence, from the third claim in the statement of Theorem~\ref{t:cm} we deduce the following. There exists a constant $\delta >0$ and a center manifold $\mathcal M^c$ such that any solution of~\eqref{e:tre}
satisfying 
\begin{equation}
\label{e:bdbd}
    |\big( U(y), p(y), \sigma \big) - \bar V| \leq \delta \quad \text{for every $y \in \R$}
\end{equation}
lies on $\mathcal M^c$. As we mentioned in Section~\ref{s:cm}, the center manifold of a system about a given equilibrium point is in general non-unique. In the following we will arbitrarily fix one center manifold and denote it by $\mathcal M^c$.

We now want to investigate the structure of $\mathcal M^c$. First, we compute its dimension, which is the dimension of the center space. By setting 
$$
    V =  \big( U, p, \sigma \big)^t \quad \text{and} \quad G(V) = \big( p,  \big[ Df (U) - \sigma I \big] p, 0 \big)^t
$$
we get that $V \in \R^{2N +1}$ and the Jacobian matrix $D G (\bar V)$ is given by
\begin{equation}
\label{e:matrice}
    D G (\bar V) =
    \left(
    \begin{array}{ccc}
               \mathbf{0} & I & \vec 0 \\
               \mathbf{0} & Df (\bar u) - \lambda_i (\bar u) I & \vec 0 \\
               0 &  0 & 0 \\
    \end{array}
    \right). \end{equation}
In the previous expression, $\mathbf{0}$ and $I$ denote the null and the identity $N \times N$ matrix respectively. 
Let us assume that the matrix $Df (\bar u)$ has $N$ real and distinct eigenvalues: this hypothesis of so-called \emph{strict hyperbolicity} is quite standard in the study of conservation laws (see for example Dafermos~\cite[Section 9.6]{Daf:book} for a discussion on the ``bad behaviors" exhibited by systems violating strict hyperbolicity). Let $\vec r_i$ denote a unit eigenvector associated to $\lambda_i(\bar u)$, namely
$$
    \big[ Df (\bar u) - \lambda_i (\bar u) \big] \vec r_i = \vec 0 \quad \text{and} \quad |\vec r_i | =1. 
$$
Then the center space of the matrix~\eqref{e:matrice} is 
$$
    V^c = \Big\{ (u,  p, \sigma) \; \text{such that} \; u \in \R^N, \; p = a \vec r_i \; \text{for some $a \in \R$, and $\sigma \in \R$}\Big\}. 
$$
Hence, the dimension of $V^c$ and of $\mathcal M^c$ is $N + 2$. This gives an answer to one of the questions we asked ourselves: if we are looking for solutions of~\eqref{e:tre} satisfying~\eqref{e:bdbd}, we need to impose $N+2$ conditions (instead of $2N +1$) to select a unique solution. 

Let us get some more precise information about the structure of $\mathcal M^c$. Let 
\begin{equation}
\label{e:delta}
    \phi_c : B(\vec 0, \delta) \cap V^c \to \R^{2N +1}
\end{equation}
be a parameterization of $\mathcal M^c$. From the proof of Theorem~\ref{t:cm} we infer that $\phi_c$ can be chosen in such a way that $\pi_c \circ \phi_c ( x_c) = x_c$ for any $x_c \in V^c$. By applying this property to system~\eqref{e:tre} we get 
\begin{equation}
\label{e:proie}
   \pi_c \circ \phi_c (U, a \vec r_i, \sigma) = (U, a \vec r_i, \sigma)
   \quad \text{for every $U \in \R^N$, $a \in \R$, $\sigma \in \R$.}  
\end{equation}
and hence 
$$
    \phi_c (U, a \vec r_i, \sigma)  = \Big(U, \psi_c (U, a, \sigma) , \sigma \Big),
$$
where $\psi_c: V^c \to \R^N$ is a suitable function whose properties we want to investigate. Let $\{ \vec r_1 \dots \vec r_N \}$ be a basis of $\R^N$ entirely made by unit eigenvectors of $Df (\bar u)$. Then we denote by $\psi_{1c} \dots \psi_{Nc}$ the corresponding components of $\psi_c$ and we get
\begin{equation}
\label{e:somma}
     \psi_c (U, a, \sigma) = \sum_{k=1}^N \psi_{kc} (U, a, \sigma) \vec r_k. 
\end{equation}
To write~\eqref{e:somma} in a more convenient form, we first exploit~\eqref{e:proie} to get 
$
   \psi_{ic} (U, a, \sigma) = a.    
$
Also, any point $(U, \vec 0, \sigma)$ is an equilibrium for~\eqref{e:tre}, so it belongs to $\mathcal M^c$. This implies that  
$$
    \vec 0 = \sum_{k=1}^N \psi_{kc} (U, a, \sigma) \vec r_k \; \textrm{and hence that}
    \quad  \psi_{kc} (U, 0, \sigma) = 0 \; \text{for every $k$, $U$ and $\sigma$}. 
$$
By exploiting the regularity of the map $\psi_{kc}$, we deduce that $\psi_{kc} (U, a, \sigma)= a \tilde \psi_{kc} (U, a, \sigma)$ for a suitable function $\tilde \psi_{kc}$. 
By combining the previous observations we then get that~\eqref{e:somma} can be rewritten as 
$$
    p = a \Big( \vec r_i + \sum_{k \neq i} \tilde \psi_{kc} (U, a, \sigma)  \vec r_k \Big).
$$
Since $\mathcal M^c$ is tangent at the equilibrium point to the center space $V^c$ we also have 
\begin{equation}
\label{e:inzero}
    \Big( \vec r_i + \sum_{k \neq i} \tilde \psi_{kc} (\bar U, 0, \lambda_i (\bar U) )  \vec r_k \Big) = \vec r_i, 
\end{equation}
namely $\tilde \psi_{kc} (\bar U, 0, \lambda_i (\bar U) ) =0$ for every $k \neq i$. 
By setting 
$$
    \hat r_i (U, a, \sigma) : =  \Big( \vec r_i + \sum_{k \neq i} \tilde \psi_{kc} (U, a, \sigma )  \vec r_k \Big) ,
$$
we then obtain that~\eqref{e:somma} can be rewritten as  
$$
     p = a \hat r_i (U, a, \sigma) . 
$$ 
For convenience, we want to reduce to the case of unit vectors, so we write the above relation as 
$$
       p = a | \hat r_i (U, a, \sigma)| \, 
       \frac{   \hat r_i (U, a, \sigma)}{|   \hat r_i (U, a, \sigma)|}, 
$$
where $|\hat r_i (U, a, \sigma)|$ denotes the length of the vector. Note that~\eqref{e:inzero} implies that $\hat r_i (\bar U, 0, \lambda_i (\bar U)) = \vec r_i$. Hence, 
the map
$$
     a \mapsto  v (a) = a | \hat r_i (U, a, \sigma)| 
$$
is locally invertible in a neighborhood of $(\bar U, 0, \lambda_i (\bar U))$: 
by assuming that the size of the constant $\delta$ in~\eqref{e:delta} is small enough, we eventually get that the manifold $\mathcal M^c$ is described by the relation
\begin{equation}
\label{e:finale}
         p = v \tilde r_i (U, v, \sigma),
\end{equation}
where we have set $\tilde r_i = \hat r_i / |\hat r_i|$. 
\subsection{Center manifold and convergence of the vanishing viscosity approximation}
\label{sus:cm:vv}
In this section we just mention the fact that center manifold techniques have been exploited by Bianchini and Bressan~\cite{BiaBrevv} to show the convergence of the solutions of the family of Cauchy problems 
\begin{equation}
\label{e:cau:eps}
    \left\{
    \begin{array}{lll}
               \ue_t + f (\ue)_x = \ee \ue_{xx} \\ 
               \ue(0, x) = u_0(x).
    \end{array} 
    \right.
\end{equation}
More precisely, the following holds. 
\begin{theorem}[Theorem 1 page 229 in~\cite{BiaBrevv}]
\label{t:biabre}
          Assume that the flux function $f: \R^N \to \R^N$ is of class $\mathcal C^2$ and that, for every $u \in \R^N$, the jacobian matrix $Df (u)$ has $N$ real and distinct eigenvalues. 
          There exists a small enough constant $\delta >0$ such that, if the initial datum in~\eqref{e:cau:eps} satisfies 
         $$
              u_0 \in L^1(\R)   \quad \text{and} \quad         \mathrm{TotVar} \{ u_0 \} \leq \delta,
         $$
         then the following properties hold.
         \begin{enumerate}
         \item For any $\ee >0$, the Cauchy problem~\eqref{e:cau:eps} admits a solution defined for any $t \ge 0$. Also, there exists a constant $C$, which does not depend neither on $\ee$ nor on $t$, such that 
         \begin{equation}
         \label{e:variazione}
             \mathrm{TotVar}_x \{ \ue (t, \cdot) \} \leq C \delta. 
         \end{equation}
         \item The solution $\ue$ is stable in $L^1$ with respect to time and to 
         $L^1$ perturbation of the initial data. Namely, if by using a semigroup 
         notation we denote 
         by $t \mapsto S^{\ee}_t u_0$ the solution of~\eqref{e:cau:eps}, then
         $$
             \| S^{\ee}_t u_0 - S^{\ee}_s v_0 \|_{L^1} \leq L \Big( \|u_0 - v_0 \|_{L^1} +
             |t -s| + |\sqrt{\ee t} - \sqrt{\ee s} | \Big) \; \text{for every $t, s \ge 0$}.
         $$
         In the previous expression, $L$ denotes a constant which does not depend neither on $\ee$ nor on the time.          
         \item As $\ee \to 0^+$, for every $t \ge 0$,  $S^{\ee}_t u_0$ converges strongly
         in $L^1_{\loc} (\R)$ to a function $S_t u_0$ providing a distributional solution of~\eqref{e:cau:eps}. Also, the semigroup $t \mapsto S_t u_0$ enjoys the stability property 
         $$
              \| S_t u_0 - S_s v_0 \|_{L^1} \leq L \Big( \|u_0 - v_0 \|_{L^1} +
             |t -s|  \Big) \; \text{for every $t, s \ge 0$}.
         $$
         
         \end{enumerate} 
         \end{theorem}
\subsection{Solutions of the Riemann problem}
\label{sus:rie}
In this section, as a direct application of the center manifold construction discussed in Section~\ref{sus:cm:come}, we go over the analysis of the solution of the so-called Riemann problem
\begin{equation}
\label{e:rie}
   \left\{
   \begin{array}{ll}
             u_t + f(u)_x =0 \\
             u(0, x) =
              \left\{
   \begin{array}{ll}
             u^+ & x> 0 \\
             u^- & x <0 \\
   \end{array}
   \right.
   \end{array}
   \right.
\end{equation}
Here $u^+$ and $u-$ are two given values in $\R^N$ and in the following we focus on the case $|u^+ - u^-| < \delta$ for a small enough constant $\delta$. The Riemann problem is a very specific example of Cauchy problem, but is nonetheless extremely important. Indeed, its study provides the building block for the construction of approximation schemes (Glimm scheme~\cite{Gli}, wave front-tracking algorithm) that have been used to establish global existence and uniqueness results for the system of conservation laws $ u_t + f(u)_x =0$: see for example Dafermos~\cite[Chapter 14]{Daf:book} and the references therein.

The first issue one has to address when studying~\eqref{e:rie} is that, in general, a distributional solution of~\eqref{e:rie} is not unique: in the attempt to select a unique solution, various {\em admissibility criteria} have been introduced. An admissibility criterium was introduced by Lax~\cite{lax}, who also constructed an \emph{admissible solution} of~\eqref{e:rie}. The analysis in~\cite{lax} relied on some technical assumptions on the flux $f$  which were later relaxed in Liu~\cite{liu:riemann}. See also Tzavaras~\cite{Tz} and Joseph and LeFloch~\cite{JosephLeFloch:olga}. Bianchini and Bressan~\cite[Section 14]{BiaBrevv} provided a general construction of an {\em admissible solution} obtained by taking the limit $\ee \to 0^+$ of the vanishing viscosity approximation
$$
     \left\{
   \begin{array}{ll}
             \ue_t + f(\ue)_x =\ee \ue_{xx} \\
             \ue(0, x) =
              \left\{
   \begin{array}{ll}
             u^+ & x> 0 \\
             u^- & x <0 \\
   \end{array}
   \right.
   \end{array}
   \right.
$$  
In both~\cite{lax, liu:riemann} and~\cite{BiaBrevv} the key point is the construction of the so-called {\em $i$-wave fan curve}. Loosely speaking, the basic idea is the following: given $u^- \in \R^N$, the $i$-wave fan curve $T_i(u^-, s)$ contains all the states $u^+$ such that an admissible solution of the Riemann problem~\eqref{e:rie} is obtained by gluing shocks (or contact discontinuities) and rarefaction waves belonging to the $i$-family. In particular, this implies that there is a self-similar, admissible solution $u(t, x) = U(x /t)$ of~\eqref{e:rie} such that $U(\xi)$ varies only in a small neighborhood of the eigenvalue $\lambda_i(u^-)$ and is constant elsewhere, $U (\xi) = u^-$ for $\xi < \lambda_i(u^-) - \delta$ and $U (\xi) = u^+$ for $\xi > \lambda_i(u^-) + \delta$.

We now informally go over the main steps of the construction in~\cite{BiaBrevv}. 
\begin{enumerate}
\item We fix $u^- \in \R^N$ and we consider the same $\tilde r_i$ as in~\eqref{e:finale}, obtained by applying the Center Manifold Theorem to system~\eqref{e:tre} about the equilibrium point $\big(u^-, \vec 0, \lambda_i(u^-) \big)$. Here $\lambda_i(u^-)$ is the $i$-eigenvalue of $Df(u^-)$.  Given $s>0$, we consider the following fixed point problem, defined on the interval $[0, s]$:  
\begin{equation}
\label{e:contra}
     \left\{
\begin{array}{lll}
      u(\tau) = u^- + {\displaystyle \int_0^{\tau} \tilde r_i (u(\xi),  v_i(\xi), 
      \sigma_i(\xi))d \xi }  \\
      v_i (\tau) = f_i (\tau, u,  v_i,  \sigma_i ) -
      \mathrm{conv} f_i (\tau )\\
      \sigma_i(\tau)=   {\displaystyle \frac{d}{d \tau}
      \mathrm{conv} f_i (\tau)}. \\
\end{array}
\right.
\end{equation}
Here 
$$
    f_i (\tau)= \int_0^{\tau}  \tilde \lambda_i (u(\xi), v_i(\xi), 
      \sigma_i(\xi))d \xi \quad \text{and} \quad \tilde \lambda_i (u(\xi),  v_i(\xi), 
      \sigma_i(\xi)) = \langle A (u) \tilde r_i,  \tilde r_i  \rangle,  
$$
namely $\tilde \lambda_i$ can be regarded as a {\em generalized eigenvalue}. Also, $\mathrm{conv} f_i (\tau)$ denotes the convex envelope of $f_i  (\tau)$, 
$$
   \mathrm{conv} f_i (\tau ) = \sup \Big\{ g (\tau) \; \text{s.t. $g$ is convex, $g(\xi) \leq f(\xi)$ for every $\xi \in [0, s]$} \Big\} .$$ By relying on the Contraction Map Theorem, one can show that~\eqref{e:contra} admits a unique continuous solution $\big(u(\tau), v_i (\tau), \sigma_i (\tau)\big)$ defined in a small enough neighborhood of $(u^-, 0, \lambda_i(u^-)$. The functions $\big(u(\tau), v_i (\tau), \sigma_i (\tau)\big)$ are all defined on the interval $[0, s]$. 
\item Assume that the interval $[a, b] \subseteq [0, s]$ satisfies the following property: for every $\tau \in [a, b]$, $f_i(\tau) =  \mathrm{conv} f_i (\tau)$. Then one can show that $\tilde r_i \big( u(\tau), v_i(\tau), 
      \sigma_i(\tau) \big) = r_i \big(u(\tau) \big)$ for every $\tau \in [a, b] $ and moreover that the solution of the Riemann problem 
$$
 \left\{
   \begin{array}{ll}
             u_t + f(u)_x =0 \\
             u(0, x) =
              \left\{
   \begin{array}{ll}
             u(b) & x> 0 \\
             u(a)& x <0 \\
   \end{array}
   \right.
   \end{array}
   \right.
$$
is a rarefaction wave of the $i$-th family (hence, in particular, the solution is continuously differentiable in $]0, + \infty[ \times \R$).  
\item
Conversely, assume that for every $\tau \in ]c, d [ \subseteq [0, s]$, $f_i(\tau) >  \mathrm{conv} f_i (\tau)$ and that $f_i(c )= \mathrm{conv} f_i (c) $, $f_i(d )= \mathrm{conv} f_i (d)$. Then, in particular, one can show that $\sigma_i$ is constant on $[c, d]$, $\sigma_i = \bar \sigma _i$. Also, $v_i (\tau) >0$ on $]c, d[$ and one can define the change of variables 
 $$
     \left\{
   \begin{array}{ll}
             \displaystyle{ \frac{dx}{dt} = \frac{1}{v_i (t)} }\\ \\
              \displaystyle{x (0) = \frac{d - c }{2}} \\
      \end{array}
      \right.         
             $$
 mapping $x: \, ]c, d[ \to ]- \infty, + \infty[$. Let $\tau (x)$ denote its inverse: one can show that the function
 $
     U(x) : = u \big( \tau (x) \big)
 $
 satisfies 
 $$
     U '' = \big( f(U) - \bar \sigma_i U \big)' \quad \text{and} 
     \quad \lim_{x \to - \infty} U(x) = u(c), \quad \lim_{x \to + \infty} U(x) = u(d). 
 $$
 In other words, there is a traveling wave connecting $u(c)$ and $u(d)$ and hence a solution of the Riemann problem
$$
 \left\{
   \begin{array}{ll}
             u_t + f(u)_x =0 \\
             u(0, x) =
              \left\{
   \begin{array}{ll}
             u(d) & x> 0 \\
             u(c)& x <0 \\
   \end{array}
   \right.
   \end{array}
   \right.
$$ 
is 
$$
    u(t, x) =  \left\{
   \begin{array}{ll}
             u(d) & x> \bar \sigma_i  t \\
             u(c)& x < \bar \sigma_i t, \\
   \end{array}
   \right.
$$
which is called either {\em shock} or {\em contact discontinuity} depending on the values of $\lambda_i \big(u(d) \big)$ and 
$\lambda_i  \big(u(c) \big)$ (see for example the book by Dafermos~\cite[Page 213]{Daf:book} for the exact definition).
\item By relying on the previous steps and on further approximation arguments, one can show that if  
$\big(u(\tau), v_i (\tau), \sigma_i (\tau)\big)$ is the solution of~\eqref{e:contra}, then the family $\ue$ solving
$$
\left\{
   \begin{array}{ll}
             \ue_t + f(\ue)_x = \ee \ue_{xx} \\
             \ue (0, x) =
              \left\{
   \begin{array}{ll}
             u^- & x <0 \\
               u(s) & x> 0 \\
   \end{array}
   \right.
   \end{array}
   \right.
$$
converges $\ee \to 0^+$ to the function
$$
    u(t, x) = 
   \left\{
   \begin{array}{ll}
               u^-             &       \mathrm{if} \; x < \sigma_i(0) t \\
               u(\tau)       &       \mathrm{if} \; x = \sigma_i(\tau) t \\
               u(s)            &      \mathrm{if} \; x > \sigma_i(s) t \\
   \end{array}
   \right.
$$
We then define the i-wave fan curve by setting $T_i (u^-, s) : = u(s)$. 
\item When $s < 0$, the construction of $T_i (u^-, s)$ is analogous, the only difference being that in~\eqref{e:contra} one has to take the concave envelope of $f_i$ instead of the convex one.
\end{enumerate}
Finally, observe that the above construction can be extended to study the solution of the Riemann problem~\eqref{e:rie} obtained by taking the limit $\ee \to 0^+$  of the more general viscous approximation
$$
      \ue_t + f(\ue)_x =\ee \Big( B (\ue) \ue_x \Big)_x,
$$
under quite general assumptions on the matrix $B$ (see Bianchini~\cite{Bia:riemann}).  For extensions to the analysis of initial-boundary value problems see Bianchini and Spinolo~\cite{BiaSpi:rie}. 
\section{ Stable Manifold Theorems}
\label{s:stable}
\subsection{The Stable Manifold Theorem}
\label{sus:stable}
Let us first recall some notations: given the ODE
\begin{equation}
\label{e:s:ode}
   {V' = G(V)}, 
\end{equation}   
let $\bar V$ be an equilibrium point, which without any loss of generality is assumed to be $\vec 0$. As before, we denote by $V^s$ and $V^u$ the stable and the unstable space respectively, defined as in~\eqref{e:stableun}. Also, as usual, $ DG (\vec 0) $ denotes the Jacobian matrix of $G$ computed at $\vec 0$. 
\begin{theorem}
\label{t:stable}
          Consider the first order ODE~\eqref{e:s:ode}
         where $G: \R^d \to \R^d$ is a $\mathcal C^2$ function satisfying 
         $G (\vec 0) = \vec 0$.  If $V^s$ is non trivial, namely $V^s \neq \{ \vec 0\}$, then there exists a constant $\delta >0$ small enough such there exists a continuously differentiable  
         stable manifold $\mathcal M^s$
         satisfying the following conditions. 
         \begin{enumerate}
         \item $\mathcal M^s$ is parameterized by a function
         $$
             \phi_s : B(\vec 0, \delta) \cap V^s \to \R^d.
         $$
       Also, $\mathcal M^s$ is tangent to $V^s$ at the origin.
         \item $\mathcal M^s$ is locally invariant for the ODE~\eqref{e:cm:ode}. Namely,  
         if $V_0 \in \mathcal M^s$, then the solution of the Cauchy problem
         \begin{equation}
         \label{e:s:cauchy}
         \left\{
         \begin{array}{lll}
                   V' = G(V) \\
                   V(0) = V_0     
         \end{array}
         \right.
         \end{equation}
         satisfies $V(t) \in \mathcal M^s$ for every $|t|$ small enough. 
         \item If $V_0  \in \mathcal M^s$ then the solution of the Cauchy problem~\eqref{e:s:cauchy} satisfies
         \begin{equation}
         \label{e:s:lim}    
             \lim_{t \to + \infty} |V(t)| e^{ct / 2 } = 0,
         \end{equation}
         where 
         \begin{equation}
         \label{e:c}
             c =  \min \{ \text{$-Re \lambda_i$ such that $\lambda_i$ is an eigenvalue of $DG (\vec 0)$, $Re \lambda_i <0$} \}
         \end{equation}
         \item Conversely, if $V(t)$ is a solution of~\eqref{e:s:ode} satisfying~\eqref{e:s:lim}, then 
         $V(t) \in \mathcal M^s$ for every $t$ large enough.   
   \end{enumerate}                  
        \end{theorem}
The proof of Theorem~\ref{t:stable} can be obtained by relying on an argument somehow similar to the one exploited in the proof of Theorem~\ref{t:cm}: a detailed proof can be found in the book by Perko~\cite{Perko} (see also the book by Katok and Hasselblatt~\cite{KatokHass} for a more general viewpoint). One should keep in mind, however, that the stable manifold (i.e., a manifold satisfying properties (1), (2), (3) and (4) in the statement of Theorem~\ref{t:stable}) is unique, while the center manifold is not. As a matter of fact, one can also define a \emph{global} stable manifold as 
$$
     \mathcal M^{gs} = \bigcup_{t \leq 0} \,
     \{ V(t) \; \text{solving~\eqref{e:s:ode} and satisfying $V(0) \in \mathcal M^s$} \}.   
$$  
In this way, we obtain a manifold which is \emph{globally invariant} for~\eqref{e:s:ode}. Namely, if $V_0 \in \mathcal M^{gs}$, then the solution $V(t)$ of~\eqref{e:s:cauchy} satisfies~$V(t) \in \mathcal M^{gs}$ for every $t$. 
   
\subsection{Applications to the analysis of the boundary layers}
\label{sus:s:app}
Let us now go back to the original problem described at the end of Section~\ref{sus:bl}: given $u_0$, we want to describe the values of $u_b$ such that problem 
$$
          \left\{
    \begin{array}{lll}
                U ' = f(U) - f (u_0) \\
                U( 0) = u_b \qquad \lim_{x \to + \infty} U(x) = u_0
    \end{array}
    \right.
$$
admits a solution. By combining parts (2) and (3) in the statement of the theorem we get that the previous problem admits a solution provided $u_b$ belongs to the stable manifold $\mathcal M^s$, whose dimension is equal to the number of the eigenvalues of $Df (u_0)$ with negative real part. If  all the eigenvalues of the Jacobian matrix $Df(u_0)$ have nonzero real part, then as a matter of fact one can also prove the converse implication: if the above problem admits a solution, then $u_b \in \mathcal M^s$. 
\subsection{The slaving manifold of a manifold of equilibria} 
\label{s:sm}
Loosely speaking, the Stable Manifold Theorem we discussed in the previous section describes all the solutions of~\eqref{e:s:ode} that when $t \to + \infty$ converge exponentially fast to the equilibrium point $\vec 0$. 

Let us consider a more general situation: assume that there exists a continuously differentiable manifold $\mathcal E \subseteq \R^d$ which contains $\vec 0$ and is entirely made by equilibria:
$$
    \forall \; V \in \mathcal E, \; G( V) = \vec 0.
$$
As a side remark, note that we are {\em not} assuming that $\mathcal E$ is {\em the} set 
of equilibria, namely we do not rule out the possibility that there are equilibria that do not belong to $\mathcal E$. In several situations, we are led to ask ourself the following question: is there any way to describe the solutions of~\eqref{e:s:ode} that converge to some point in $\mathcal E$? Note that by applying the Stable Manifold Theorem we describe the solutions that converge to a given point, here we just require that the limit exists and that it lies somewhere on $\mathcal E$. 

A positive answer to the previous question is provided by the theorem below.
\begin{theorem}
\label{t:us}
         Consider the first order ODE~\eqref{e:s:ode}
         where $G: \R^d \to \R^d$ is a $\mathcal C^2$ function satisfying 
         $G (\vec 0) = \vec 0$. Denote by $E$ the space tangent to $\mathcal E$ at the origin. Also, as before let $V^s$ be the stable space of $G$ at the origin, defined as in~\eqref{e:stableun}: assume that $V^s$ is non trivial, namely $V^s \neq \{ \vec 0\}$. 
     Finally, let the constant $c$ be as in~\eqref{e:c}.
Then there exists a constant $\delta >0$ small enough and a continuously differentiable  
         manifold $\mathcal M^{us}$ 
         satisfying the following conditions.
         \begin{enumerate}
         \item $\mathcal M^{us}$ is parameterized by a function
         $$
             \phi_{us} : B(\vec 0, \delta) \cap \Big( V^s  \oplus E \Big)    \to \R^d.
         $$
         and it is tangent to $V^s \oplus E$ at the origin.
         \item $\mathcal M^{us}$ is locally invariant for the ODE~\eqref{e:cm:ode}. Namely,  
         if $U_0 \in \mathcal M^{us}$, then the solution of the Cauchy problem
         \begin{equation}
         \label{e:sl:cauchy}
         \left\{
         \begin{array}{lll}
                   V' = G(V) \\
                   V(0) = V_0     
         \end{array}
         \right.
         \end{equation}
         satisfies $V(t) \in \mathcal M^{us}$ for every $|t|$ small enough. 
         \item If $V_0  \in \mathcal M^{us}$ then there exists 
         $V_{\infty} \in \mathcal E$ such that the  solution of the 
         Cauchy problem~\eqref{e:sl:cauchy} satisfies
         \begin{equation}
         \label{e:sl:lim}    
             \lim_{t \to + \infty} |V(t) - V_{\infty} | e^{ c t / 2} = 0.
         \end{equation}
         \item Conversely, if $V(t)$ is a solution of~\eqref{e:s:ode} 
         satisfying~\eqref{e:sl:lim} for some $|V_{\infty}| \leq \delta $, then 
         $V(t) \in \mathcal M^{us}$ for every $t$ large enough.   
   \end{enumerate}                  
        \end{theorem}
 The manifold $\mathcal M^{us}$ can be regarded as a {\em uniformly stable manifold} because it contains the solutions of~\eqref{e:s:ode} that converge to some point in $\mathcal E$ with exponential speed uniformly bounded from below. 
As a matter of fact, the idea behind the definition of uniformly stable manifold is a particular case of the idea of {\em slaving manifold}. Loosely speaking, the slaving manifold relative to some manifold $\mathcal A$ contains all the solutions of~\eqref{e:s:ode} that when $t \to +\infty$ approach $\mathcal A$ with exponential speed: see again the book by Katok and Hasselblatt~\cite{KatokHass}. We came back to this in Section~\ref{s:sode} by considering in Lemma~\ref{l:slaving} another class of slaving manifolds.        

Theorem~\ref{t:us} can be proven by relying on the Hadamard Perron Theorem: see Ancona and Bianchini~\cite{AnBiapro} and the book by Katok and Hasselblatt~\cite{KatokHass}. 
\section{ Invariant manifolds for a class of singular ordinary differential equations}          
\label{s:sode}
This section focuses on a class of ordinary differential equations in the form
\begin{equation}
\label{e:sode}
     \frac{d V}{d t} = \frac{1}{\zeta (V)} F(V),
\end{equation}
where $V \in \R^d$, $F: \R^d \to \R^d$ is a smooth function and $\zeta$ is a real valued, smooth function. The singularity of the equation comes from the fact that $\zeta$ can attain the value $0$. 
\subsection{Motivations}
\label{sus:sode:mot}
Singular ODEs like~\eqref{e:sode} may arise in the analysis of the viscous approximation 
$$
     u_t + f(u)_x = \big( B(u) u_{x} \big)_x
$$
when the $N \times N$ matrix $B(u)$ does not have full rank, which is the case of most of the physically relevant examples. 
In particular, let us focus on the case of the compressible Navier-Stokes equation in one space variable:
\begin{equation}
\label{e:ns:eul}
       \left\{
       \begin{array}{lll}
              \rho_t + ( \rho v )_x =0 \\
	      (\rho v)_t + \Big( \rho v^2 + p \Big)_x = \displaystyle{ \Big( \nu v_x  \Big)_x } \\
	      \displaystyle{ \Big( \rho e + \rho \frac{v^2}{2}\Big)_t + \Big(v \Big[ \frac{1}{2} \rho v^2 
	      + \rho e + p \Big] \Big)_x = \Big( k \theta_x + 
	      \nu v v_x \Big)_x}. \\
       \end{array}
       \right.
\end{equation}
Here the unknowns $\rho(t, \, x), \, v(t, \, x)$ and $\theta(t, \, x)$ are the density of the fluid, the velocity of the particles in the fluid and the absolute temperature respectively. The function $p= p(\rho, \, \theta) >0$ is the pressure and satisfies $p_{\rho} >0$, while $e$ is the internal energy. In the following we will focus on the case of a polytropic gas, so that $e$ satisfies $
   e = R   \theta /( \gamma -1 )$,  
where $R$ is the universal gas constant and $\gamma >1$ is a constant specific of the gas. Finally, $\nu(\rho)>0$ and $k(\rho)>0$ denote the viscosity and the heat conduction coefficients respectively. 

Let us now compute the equations satisfied by the traveling waves and the steady solutions of~\eqref{e:ns:eul}: for simplicity, in the following we focus on steady solutions only, but the considerations we make can be repeated with small changes in the case of traveling waves.  We have to consider the equation
\begin{equation}
\label{e:ns:eul:st}
       \left\{
       \begin{array}{lll}
               ( \rho v )_x =0 \\
	       \Big( \rho v^2 + p \Big)_x = \displaystyle{ \Big( \nu v_x  \Big)_x } \\
	      \displaystyle{  \Big(v \Big[ \frac{1}{2} \rho v^2 
	      + \rho e + p \Big] \Big)_x = \Big( k \theta_x + 
	      \nu v v_x \Big)_x}. \\
       \end{array}
       \right.
\end{equation}
To write~\eqref{e:ns:eul:st} in a convenient form we set 
$$
    \vec w = ( v_x, \, \theta_x)  
$$
and hence we get
\begin{equation}
\label{e:block}
     \left(
    \begin{array}{ccc}
              v             &    A_{12} \\
              A_{21}    &    A_{22} \\
    \end{array}
    \right)
    \,
     \left(
    \begin{array}{ccc}
             \rho_x  \\
              \vec w \\
    \end{array}
    \right)  = 
    \,
     \left(
    \begin{array}{ccc}
              0            &    \vec 0\\
              \vec 0   &    b  \\
    \end{array}
    \right)
     \left(
    \begin{array}{ccc}
             \rho_{xx}  \\
              \vec w_x \\
    \end{array}
    \right) 
\end{equation}
Here, $A_{12}  \in \R^2$ is a row vector, $ A_{21} \in 
\R^2 $ is a column vector  and they both depend on $(\rho, v, \theta, \vec w)$. Finally, $A_{22}$ and $b$ are both $2 \times 2$ matrices depending on $(\rho, v, \theta, \vec w)$. The exact expression of these terms is not important here: the important thing to keep in mind is that the block $b$ is invertible. 

To write~\eqref{e:block} in a more explicit form, we first assume $v\neq 0$ and from the first line we get 
$$
     \rho_x = - \frac{1}{v} \, A_{21} \cdot \vec w,
$$ 
where $\cdot$ denotes the standard scalar product. By plugging this expression in the last two lines of~\eqref{e:block} we get 
$$
    z_x = b^{-1} \left[  A_{22} - \frac{1}{v} \, A_{12} A_{21}   \right] \vec w.
$$
 Note that this expression makes sense since $b$ is invertible. By combining the previous considerations we get that~\eqref{e:ns:eul:st} can be written in the form~\eqref{e:sode} provided that 
 $$
 V = 
 \left(
 \begin{array}{ccc}
           \rho \\
            v \\
            \theta \\
            \vec w \\
 \end{array}
 \right) \in \R^5, \qquad
 F(V) = 
 \left(
 \begin{array}{ccc}
        - A_{12} \cdot \vec w \\
         v  \, \vec w\\
            b^{-1} \left[  A_{22} v  - A_{21} A_{12}   \right] \vec w
           \end{array}
 \right) \in \R^5 
 $$
 and $ \zeta (V) = v \in \R$. Note that $v$ is the velocity of the fluid and in general it can attain the value $0$. 
 \begin{remark}
 Equation~\eqref{e:ns:eul} is the compressible Navier-Stokes equation written in Eulerian coordinates. If we write the Navier-Stokes equation using Lagrangian coordinates and we compute the equation satisfied by the traveling waves and by the steady solutions, we obtain standard ODEs with no singularity. See for example Rousset~\cite{Rousset:char} for related  analysis. 
 \end{remark}
 \subsection{What in principle can go wrong with~\eqref{e:sode}}
 \label{s:counterxamples}
 There is a wide literature concerning the analysis of equations in the form 
 \begin{equation}
 \label{e:sodepar}
     \frac{d V}{ dt} = \frac{1}{\ee} \, F^{\ee} (V), \quad V \in \R^d,
 \end{equation}
 where $\ee$ is a parameter that can go to $0$ and $F^{\ee}$ is a family of smooth functions. For an overview, see for example the lecture notes by Jones~\cite{Jones}, which in particular provide an overview of Fenichel's works~\cite{Fenichel_asy, Fenichel}. 
 The reason why we cannot directly apply these results to the analysis of~\eqref{e:sode} is because the singularity $\zeta (V)$ depends on the solution $V$ itself, while in~\eqref{e:sodepar} the singularity is just a parameter. Hence, when handling~\eqref{e:sode} we have to tackle the possibility that $z(V)$ is nonzero at time $t=0$, but then $z (V)$ reaches the singular value $0$ in finite time. When this happens, the solution may experience a loss of regularity, as the example below shows. 
 \begin{example}
 \label{ex:noreg}
 Let us consider the ODE
  \begin{equation}
           \label{e:ex:fast}
           \left\{
           \begin{array}{ll}
                      d v_1 / dt = - v_2 / v_1 \\
                      d v_2 / dt = - v_2 . \\
           \end{array}
           \right.                     
 \end{equation}    
It can be written in form \eqref{e:sode} provided that $V= (v_1, \, v_2)^t$, $\zeta (V) = v_1$ and 
$$
    F (V) =
    \left(
    \begin{array}{cc}
                - v_2 \\
                -v_2 v_1  \\
     \end{array}
     \right)
$$
The solution of \eqref{e:ex:fast} is 
\begin{equation}
\label{e:ex:sol}
     \left\{
    \begin{array}{lll}
                \displaystyle{ v_1 (t) = \sqrt{  v_1 (0) +  v_2 (0) \big(  e^{- t } - 1 \big)  }} \\
                \\
                   \displaystyle{  v_2 (t) = v_2 (0) e^{- t } } \\
    \end{array}
    \right.
\end{equation}
By choosing $v_2 (0) > v_1 (0) >0$, we get that $\zeta (V) =v_1 (t)$ can attain the singular value $0$ for a finite $t$. Note that at that point $t$ the first derivative $d v_1 / dt$
blows up: thus, in particular, the solution \eqref{e:ex:sol} of \eqref{e:ex:fast} is not $\mathcal C^1$.  
\end{example}
\subsection{Goals in studying~\eqref{e:sode}}
\label{s:goal}
The analysis concerning~\eqref{e:sode} and discussed in the following aims at applying the ODE techniques described in the previous sections to the study of the viscous profiles of the Navier-Stokes equation.  

More precisely, the goals are the following.
\begin{enumerate}
\item Find possible extensions of the Center Manifold Theorem~\ref{t:cm} and of the Uniformly Stable Manifold Theorem~\ref{t:us} to the case of the singular ODE~\eqref{e:sode}. Note that both theorems provide \emph{local} results, so to study this extension we can restrict to the solutions of~\eqref{e:sode} belonging to a small enough neighborhood of an equilibrium point.
\item Rule out the losses of regularity discussed in Example~\ref{ex:noreg}. To do so, we have to find conditions that guarantee that the following property is satisfied.
\begin{equation}
\label{e:prop}
         \text{ If $\zeta (V) \neq 0$ at $t=0$, then $\zeta (V) \neq 0$ for every $t \ge 0$}.
\end{equation}  
\item Impose on~\eqref{e:sode} hypotheses satisfied by the viscous profiles of the Navier-Stokes equation. 
\end{enumerate}

\subsection{Hypotheses}
\label{sus:hyp} 
This section lists the conditions imposed on the functions $F$ and $\zeta$ to study~\eqref{e:sode}, and comments on them.
\begin{hyp}
\label{h:reg}
         The functions $F: \R^d \to \R^d$ and $\zeta: \R^d \to \R$ are both regular (being of class $\mathcal C^2$ is enough). Also, $F( \vec 0) = \vec 0$ and $\zeta (\vec 0) = 0$. 
\end{hyp}
The second part of Hypothesis~\ref{h:reg} simply says that we are assuming that our equation is singular and that the equilibrium point is $\vec 0$, which of course is not restrictive.

\begin{hyp}
\label{h:zeta}
          We have $\nabla \zeta (\vec 0) \neq \vec 0$. 
\end{hyp}
Let $\mathcal S$ be the singular set
\begin{equation*}
         \mathcal S : = \big\{ V: \; \zeta (V) =0   \big\}.
\end{equation*}
By relying on the Implicit Function Theorem and on Hypothesis~\ref{h:zeta} we get that in a small enough 
neighbourhood of $\vec 0$ the set $\mathcal S$ is actually a manifold of dimension $d-1$, where $d$ is the dimension of $V$.  
\begin{hyp}
\label{h:center}
          Let $\mathcal M^c$ be any center manifold for 
          \begin{equation}
           \label{e:nsin2}
                    \frac{d V}{ d \tau} = F(V)
          \end{equation}
          about the equilibrium point $\vec 0$. If $|V|$ is sufficiently small and $V$ belongs to the intersection $\mathcal M^c \cap \mathcal S$ , then $V$ is an equilibrium for~\eqref{e:nsin2}, namely $F(V)= \vec 0$ . 
\end{hyp}          
To understand why we impose Hypothesis~\ref{h:center} let us consider the following simple linear example.
\begin{example}
 \begin{equation}
 \label{e:giro}
       \left\{
           \begin{array}{ll}
                      d v_1 / d \tau =  v_2 / \ee  \\
                      d v_2 / d \tau = - v_1 / \ee   \\
                      d \ee / d \tau = 0. 
           \end{array}
           \right.   
 \end{equation}
 The first component of the solution is 
 $$
     v_1 (t) = A \cos (t / \ee) + B \sin (t / \ee), 
 $$
 where $A$ and $B$ are real parameters. Letting $\ee \to 0^+$, we get that in general if $t \neq 0$ there exists no pointwise limit of $v_1$. Note that~\eqref{e:giro} does not satisfy Hypothesis~\ref{h:center}. Indeed, the center space is the whole $\mathbb R^3$ and it coincides with the center manifold. However, 
 $$
     F(V) =   
     \left(
           \begin{array}{ccc}
                        v_2  \\
                     - v_1    \\
                     0
           \end{array}
           \right)   
 $$ 
 is not identically zero when $\ee =0$.         
 \end{example}
 We now assume that there exists a non-degenerate set of equilibria.   
 \begin{hyp}
 \label{h:tras}
         There exists  a one-dimensional manifold of equilibria  $\mathcal E$ for the equation
         $$
              \frac{d V}{ d \tau} = F(V)
         $$
         which contains $\vec 0$ and which is transversal to $\mathcal S$.
\end{hyp}
The one-dimensional manifold $\mathcal E$ is transversal to the manifold $\mathcal S$ if the intersection $\mathcal S \cap \mathcal M^{eq}$ has dimension zero. 
\begin{hyp}
\label{h:fast}
          For every $V \in \mathcal S$,  
          \begin{equation*}
                    \nabla \zeta (V) \cdot F(V) =0  .   
          \end{equation*}
\end{hyp}
Hypothesis~\ref{h:fast} is needed to rule out losses of regularity like the one in \eqref{e:ex:fast}. Indeed, let us go back to example \eqref{e:ex:fast}:  we have
$$
    F (v_1, v_2) =  \left(
           \begin{array}{ll}
                      - v_2  \\
                   - v_2 v_1 \\
           \end{array}
           \right), \quad \zeta (v_1, v_2)= v_1  \quad  \text{and} \quad \mathcal S = \{ (v_1, v_2): \; v_1 =0  \},
$$
hence 
$$
    \nabla \zeta (v_1, v_2) \cdot F (v_1, v_2) = - v_2,  
    $$
which in general is nonzero on $\mathcal S$. 

To introduce the last hypothesis, let us consider the function 
$$
    G(V) := \frac{   \nabla \zeta (V) \cdot F(V)}{ \zeta (V)}.
$$
Thanks to Hypothesis~\ref{h:fast} and to the regularity of the functions $\zeta$ and $F$, the function $G$
can be defined by continuity on the surface $\mathcal S$. 
\begin{hyp}
\label{h:slow}
          Let $U \in \mathcal S$ be an equilibrium for \eqref{e:nsin2}, namely $\zeta (V)= 0$ and $F(V) = \vec 0$. Then 
          \begin{equation*}
                            G(V) = 0. 
          \end{equation*}
\end{hyp}
Like Hypothesis~\ref{h:fast}, Hypothesis~\ref{h:slow} is needed to rule out losses of derivatives.  Indeed, let us consider the following example.
\begin{example}
Let us consider the ODE
           \begin{equation}
           \label{e:ex:slow}
           \left\{
           \begin{array}{lll}
                      d v_1 / dt = - v_3 \\
                      d v_2 / dt = - v_2 / v_1  \\
                      d v_3 / dt = - v_3, \\
           \end{array}
           \right.           
           \end{equation}         
then we can set 
$$
    F( v_1, \, v_2, \, v_3 ) =  \left(
    \begin{array}{ccc}
                  - v_3 v_1   \\
                - v_2  \\
                 - v_3 v_1 \\
     \end{array}
     \right),
      \; \zeta ( v_1, \, v_2, \, v_3 ) = v_1 \;  \text{and} \; \mathcal S = \{ (v_1, v_2, v_3): \; v_1 =0  \}, 
$$
so that 
$$
   \nabla \zeta \cdot F = - v_3 v_1 =0 \quad \text{on $\mathcal S$}   
$$
and hence Hypothesis~\ref{h:fast} is satisfied. As a matter of fact, one can verify that Hypotheses~\ref{h:reg}~$\dots$~\ref{h:tras} are all satisfied as well. However, Hypothesis~\ref{h:slow} is violated: indeed, 
$$
  G ( v_1, \, v_2, \, v_3 ) = - v_3,
$$
which is in general different from $0$ on 
$$
     \big\{ ( v_1, \, v_2, \, v_3 ): \text{such that $F  ( v_1, \, v_2, \, v_3 ) = \vec 0$ and 
     $\zeta  ( v_1, \, v_2, \, v_3 )  = 0$} \big\}   =
     \big\{ (0, 0, v_3) \big\} . 
$$
By solving~\eqref{e:ex:slow} explicitly we get, after some computations, 
$$
\left\{
\begin{array}{ll}
         &  \displaystyle{ v_1 (t) = v_1 (0) (1- A ) + A v_1 (0) e^{-t}  }  \\
         &  \displaystyle{  v_2 (t) = B \Big[(1-A) e^t + A \Big]^{ \big[  (A-1) v_1 (0) \big]^{-1}} } \\
         &   \displaystyle{ v_3 (t) = A v_1 (0) e^{- t }},
          \phantom{ \displaystyle{  v_2 (t) = B \Big[(1-A) e^t + A \Big]^{ \frac{1 }{  (A-1) v_1 (0) }} } }  \\
\end{array}
\right.
$$       
where $A$ and $B$ are suitable constants depending on the initial data.
If $(A-1)v_1 (0) >1$, then the first derivative $d v_2 / dt$ blows up at $t = \ln (A / A-1)$. Note that this is exactly the value of $t$ at which $v_1 (t)$ attains $0$.  
In general, for every $v_1 (0) > 0$ if $1/ (A-1)v_1 (0) $ is not a natural number, then the solution is not in $\mathcal{C}^m $ for $m= [1 / (A-1)v_1 (0)  ] +1 $. Here $ [1 / (A-1)v_1 (0) ] $ denotes the entire part.  Thus,  we have a loss of regularity in higher derivatives.
\end{example} 
\begin{remark}
\label{r:ns}
        One can verify that Hypotheses~\ref{h:reg}$\dots$~\ref{h:slow} are all verified by the viscous profiles (travelling waves and steady solutions) of the Navier Stokes equation~\eqref{e:ns:eul}. 
\end{remark}

\subsection{A toy model}
\label{s:toy}
To get a flavor of the behaviors we expect to encounter, we start by considering the case when the singularity is only a parameter. Actually, to simplify things we focus on a toy model where everything is linear and we can carry out the computations explicitly. Let us consider the system 
\begin{equation}
\label{e:toy}
    \left\{
    \begin{array}{ll}
                      d v_1 / dt = -  5 v_1  \\
                      d v_2 / dt = - v_2 / \ee \\
                      d \ee / dt = 0.          
    \end{array}
    \right. 
\end{equation}
First, we note that the subspace 
$$ 
   E = \{ (0, \, 0, \, \ee ): \; \ee \in \mathbb R \}
$$ is entirely made by equilibria. If $\ee >0$ the uniformly stable manifold relative to $E$ is the whole space $\mathbb R^3$, since any solution of \eqref{e:toy} decays exponentially fast to a point in $E$. However,  
$$
    v_1 \sim e^{- 5 t} \quad \text{and} \quad v_2 \sim e^{-t / \ee}.
$$
In other words, the speed of exponential decay of $v_1$ is bounded with respect to $\ee$ and $v_1$ is not affected by the presence of the singularity. Conversely, the speed of exponential decay of $v_2$ gets faster and faster as $\ee \to 0^+$. The second component of the solution is thus strongly affected by the presence of the singularity and can be regarded as a \emph{fast dynamic}, while the first component is a \emph{slow dynamic}.  

Summing up, in~\eqref{e:toy} any orbit lying on the uniformly stable manifold relative to $E$ decomposes as the sum of a slow and a fast dynamic. 

To see what happens if we look for extensions of the Center Manifold Theorem, let us modify~\eqref{e:toy} by adding two more equations:
\begin{equation}
\label{e:toy2}
 \left\{
    \begin{array}{ll}
                      d v_1 / dt = -  5 v_1  \\
                      d v_2 / dt = - v_2 / \ee \\
                      d v_3 / dt = - v_4  \\
                      d v_4 / dt = v_3 \\
                      d \ee / dt = 0.          
    \end{array}
    \right. 
\end{equation}
In this way, we have that 
$$
    v_3 (t), v_4(t)  \sim A \cos t  + B \sin t 
$$
and hence we have a nontrivial center manifold. Note that we could not have added something like 
$$
\left\{
\begin{array}{ll}
        d v_3 / dt = - v_4 / \ee  \\
                      d v_4 / dt = v_3 / \ee \\
\end{array}
\right.
$$
because this would have violated Hypothesis~\ref{h:center}.  In other words, the center manifold can contain only {\em slow dynamics}.

\subsection{Slow and fast dynamics in the general case}
\label{sus:slow}
Singling out the slow and the fast dynamics and defining the center and the uniformly stable manifold for the toy model~\eqref{e:toy2} is a trivial task. However, to understand how we can proceed in the general case, let us first make some preliminary considerations, still focusing on~\eqref{e:toy2}.
\begin{enumerate}
\item We introduce the change of variables
$$
    \tau = t / \ee.
$$
Then~\eqref{e:toy2} becomes 
\begin{equation}
\label{e:toytau}
  \left\{
    \begin{array}{ll}
                      d v_1 / d \tau= -  5 v_1 \ee  \\
                      d v_2 / d \tau = - v_2 \\
                      d v_3 / d \tau = - v_4 \ee \\
                      d v_4 / d \tau = v_3  \ee \\
                      d \ee / d \tau = 0.          
    \end{array}
    \right. 
\end{equation}
 \item We can single out the {\em slow dynamics} of~\eqref{e:toy2} by considering the {\em center space} of~\eqref{e:toytau}, which is given by $\{ (v_1, 0, v_3, v_4, \ee) \}$.
\item Let us now restrict the analysis to the center space of~\eqref{e:toytau}: we obtain
$$
 \left\{
    \begin{array}{ll}
                      d v_1 / d \tau= -  5 v_1 \ee  \\
                      d v_3 / d \tau = - v_4 \ee \\
                      d v_4 / d \tau = v_3  \ee \\
                      d \ee / d \tau = 0.          
    \end{array}
    \right. 
$$ 
Hence, we can come back to the original variable $t$, obtaining an equation with no singularity in it:
\begin{equation}
\label{e:slow}
 \left\{
    \begin{array}{ll}
                      d v_1 / d t= -  5 v_1   \\
                      d v_3 / d t = - v_4  \\
                      d v_4 / d t = v_3   \\
                      d \ee / d t = 0.          
    \end{array}
    \right. 
\end{equation}
\item We can then consider the {\em center space} of~\eqref{e:slow}, which is given by 
$V^{c} = \{ (0, v_3, v_4, \ee) \}$. Then the center manifold of the original system~\eqref{e:toy2} is exactly $V^c$. Also, we can observe that $E = \{(0, 0, 0, \ee) \}$ is a subspace of equilibria for~\eqref{e:slow} and we can consider the uniformly stable space of~\eqref{e:slow} with respect to $E$: it is given by $\{(v_1, 0, 0, \ee) \}$.  
\item Let us now consider the uniformly stable manifold of~\eqref{e:toy2} with respect to the manifold of equilibria $E = \{(0, 0, 0, \ee) \}$. It is given by 
$
    \{ (v_1, v_2, 0, 0, \ee) \}: 
$  
as mentioned before, any solution in there writes as a fast dynamic (in this case, $(0, v_2, 0, 0)$) plus a slow dynamic (in this case, $(v_1, 0, 0, \ee)$) belonging to the uniformly stable manifold of~\eqref{e:slow} with respect to $E$. 
\end{enumerate}
Let us now discuss how one can proceed in the general case 
\begin{equation}
\label{e:sode3}
    \frac{d V}{ d t } = \frac{1}{\zeta (V)} F(V). 
\end{equation}
Actually, we will not enter into the technical details, which can be found in~\cite{BiaSpi:ode}. 
The heuristic idea is proceeding by following the same steps (1)$\dots$(5) as in the case of the toy model. Of course, now we have to handle additional difficulties coming from the fact that there is a non linearity in the function $F(V)$. Also, $\zeta (V)$ is no more a parameter, but depends on the solution 
itself, and hence we have to rule out the possibility that $\zeta (V) \neq 0$ at $t=0$ but $\zeta (V)$ reaches $0$ in finite time. 

Actually, as a preliminary step it is convenient to introduce a change of coordinates which allows to write~\eqref{e:sode3} in a more convenient form, however, this is a mostly technical step and in here we do not go into the details (they can be found in~\cite[Proposition 4.1]{BiaSpi:ode}).

Let us now informally discuss how to extend the previous steps to the general case~\eqref{e:sode3}.  
\begin{enumerate}
\item We would like to use the variable $\tau$ defined by solving the Cauchy problem
\begin{equation}
\label{e:cau}
\left\{
\begin{array}{ll}
           d \tau / dt = \zeta (V(t)) \\
           \tau (0) = 0. \\
\end{array}
\right.
\end{equation}
However,~\eqref{e:cau} defines a change of variables 
$\tau: ]- \infty;  + \infty [ \to ]- \infty;  + \infty [ $ provided that $\zeta (V(t))$ does not attain the value $0$. Hence, we actually proceed as follows: we study the equation
\begin{equation}
\label{e:tau}
           \frac{d V}{ d \tau } =  F(V)
\end{equation} 
which is {\em formally} obtained from~\eqref{e:sode3} by relying on~\eqref{e:cau}. In particular, we focus on the solutions of~\eqref{e:tau} lying on suitable invariant manifolds. We then prove that {\em for these solutions}, the Cauchy problem 
\begin{equation}
\label{e:cau2}
   \left\{
\begin{array}{ll}
           d t / d \tau= 1 / \zeta (V(\tau )) \\
           t (0) = 0. \\
\end{array}
\right.
\end{equation}
 defines a smooth diffeomorphism $t :  ]- \infty;  + \infty [ \to ]- \infty;  + \infty [ $. Hence, we can come back to the original equation~\eqref{e:sode3} and {\em a posteriori} infer that the change of variables~\eqref{e:cau} is well justified and it is inverse of the one defined by~\eqref{e:cau2}, so our argument works.  
 \item Hypothesis~\ref{h:reg} implies that $F (\vec 0) = \vec 0$. We then consider the ODE~\eqref{e:tau} linearized about the equilibrium point $\vec 0$.  
We then say that a {\em manifold of the slow dynamics} is a center manifold of~\eqref{e:toy2} about the equilibrium point $\vec 0$ (any center manifold works).   
 \item By exploiting Hypothesis~\ref{h:slow}, one can show that if $V(\tau)$ is a solution of~\eqref{e:tau} lying on the manifold of the slow dynamics, then the Cauchy problem~\eqref{e:cau2} defines a continuously differentiable diffeomorphism. Hence by restricting to the manifold of the slow dynamics we get that~\eqref{e:sode3} and~\eqref{e:tau} are equivalent. In particular, by relying on Hypothesis~\ref{h:center}  we can show that the ODE satisfied by the solutions $V(t)$ lying on the manifold of the slow dynamics is actually nonsingular.  
\item  We then consider system~\eqref{e:sode3} restricted on the manifold of the slow dynamics and we consider a center manifold about the equilibrium point $\vec 0$: we can directly apply Theorem~\ref{t:cm} because this system is nonsingular. This manifold can be regarded as a center manifold of the original system~\eqref{e:sode3}: more precisely, we get the following theorem.
\begin{theorem}
\label{t:sode:center}
         Assume that Hypotheses~\ref{h:reg} $\dots$~\ref{h:slow} are all satisfied. Then there exists a manifold $\mathcal M^c$ defined in a small enough neighborhood of $\vec 0$ and satisfying the following conditions:
         \begin{enumerate}
         \item if $U(t)$ is a solution of~\eqref{e:sode3} lying on $\mathcal M^c$ and $\zeta (V) \neq 0 $ at $t =0$ then $\zeta (V) \neq 0$ for every $t \in \R$.  
         \item $\mathcal M^c$ s locally invariant for~\eqref{e:sode3}. Namely, if $\bar V \in \mathcal M^c$, then the solution of the Cauchy problem 
         $$
         \left\{
         \begin{array}{ll}
                    \displaystyle{ \frac{d V}{ d t} = \frac{1}{ \zeta (V)} F(V)     } \\ \\
                    V(0) = \bar V \\
         \end{array}
         \right. 
         $$
         satisfies $V(t) \in \mathcal M^c$ for $|t|$ small enough.
         \item There exists a small enough constant $\delta >0$ such that the following holds. If $V(t)$ is a solution of~\eqref{e:sode3} satisfying $| V(t)| \leq \delta $ for every $t \in \R$, then $V(t) \in \mathcal M^c$. 
         \end{enumerate}  
\end{theorem}
\begin{remark}
In the statement of the Center Manifold Theorem~\ref{t:cm} and of the Uniformly Stable Manifold Theorem~\ref{t:us} one assumes that suitable subspaces (the center and the stable subspace, respectively) are non trivial. One should actually add an analogous condition to the statement of Theorem~\ref{t:sode:center}. However, to simplify the notations we can use the convention that $[F(\vec 0) / \zeta (\vec 0)] = \vec 0 / 0 : = \vec 0$ and that if a subspace is trivial then any manifold parameterized by this subspace is just the origin $\vec 0$. With these conventions the statement of Theorem~\ref{t:sode:center} is formally correct. We will adopt the same convention in the statements of Theorem~\ref{t:sode:stable}. \end{remark}
 \item  We now want find an extension of the Uniformly Stable Manifold Theorem. First, we recall that Hypothesis~\ref{h:tras} 
ensures that there is a manifold of equilibria $\mathcal E$ for~\eqref{e:tau} which is transversal to the hypersurface $\{ V: \zeta (V) =0\}$. One can show that  the points of $\mathcal E$ are also equilibria for system~\eqref{e:sode3} restricted on the manifold of the slow dynamics. 

 By relying on the analogy with the toy model, we would like to define a manifold $\mathcal M^{us}$ such that any solution lying on $\mathcal M^{us}$ decays exponentially fast to some point in $\mathcal E$. Also, any orbit lying on $\mathcal M^{us}$ should be written as the sum of a fast dynamic and of some stable component coming from the slow dynamics. However, since we are in the non linear case we expect some interactions between the different components. The precise result is the following. 
\begin{theorem}
\label{t:sode:stable}
         Let Hypotheses~\ref{h:reg}$\dots$~\ref{h:slow} hold and denote by $\mathcal E$ the same manifold of equilibria as in Hypothesis~\ref{h:tras}. Then there exist a constant $c>0$ and a manifold $\mathcal M^{us}$ which is defined in 
         a small enough neighborhood of $\vec 0$ and satisfies the following conditions.
         \begin{enumerate}
         \item If $V(t)$ is a solution of~\eqref{e:sode3} lying on $\mathcal M^{us}$ and $\zeta (V) \neq 0 $ at $t =0$ then $\zeta (V) \neq 0$ for every $t \ge 0.$ 
         \item $\mathcal M^{us}$ is locally invariant for~\eqref{e:sode3}. Namely, if $\bar V \in \mathcal M^{us}$, then the solution of the Cauchy problem 
         $$
         \left\{
         \begin{array}{ll}
                    \displaystyle{ \frac{d V}{ d t} = \frac{1}{ \zeta (V)} F(V)     } \\ \\
                    V(0) = \bar V \\
         \end{array}
         \right. 
         $$
         satisfies $V(t) \in \mathcal M^{us}$ for $|t|$ small enough. 
         \item If $V(t)$ lies on $\mathcal M^{us}$ then there exists $V_{\infty} \in \mathcal E$ such that
         $$
             \lim_{t \to + \infty} |V(t) - V_{\infty}| e^{c t} =0. 
               $$
      \item If $V(t)$   
      lies on $\mathcal M^{us}$ then
      $$
            V(t) = V_{sl}(t) + V_f (t) + V_p(t),
      $$
      where $V_{sl}$ lies on the manifold of the slow dynamics and $V_f$ satisfies 
      $$
         \lim_{\tau \to + \infty} | V_f (\tau) |e^{ c \tau} =0.
      $$
      Here the variable $\tau$ is defined by solving the Cauchy problem~\eqref{e:cau2}. 
      Finally, $V_p$ is an interaction term which is small with respect to the previous two (in the sense specified in the statement of~\cite[Theorem 4.2]{BiaSpi:ode}).  
               \end{enumerate}  
\end{theorem}
The proof of Theorem~\ref{t:sode:stable} relies on the following lemma of standard flavor.  
\begin{lemma}
\label{l:slaving} 
         Consider the ODE
         \begin{equation}
         \label{e:sode:f}
             \frac{dV}{d\tau} = F(V)
         \end{equation}
         and assume that Hypothesis~\ref{h:reg} is satisfied and that the center and the stable space about the equilibrium point $\vec 0$ are both non-trivial . Let $\mathcal S_0$ be a submanifold of a center manifold of~\eqref{e:sode:f} about the equilibrium point $\vec 0$. Also, assume that $\mathcal S_0$ is locally invariant for~\eqref{e:sode:f}. Then there exists a constant $c >0$ and a manifold $\mathcal M$ which is defined in a small enough neighborhood of $\vec 0$, is locally invariant for~\eqref{e:sode:f} and satisfies the following properties. If $V(\tau)$   
      lies on $\mathcal M$ then
      $$
            V(\tau) = V_0(\tau) + V_f (\tau) + V_p(\tau),
      $$
      where $V_0$ lies on $\mathcal S_0$ and $V_f$ satisfies 
      $$
         \lim_{\tau \to + \infty} | V_f (\tau) |e^{ c \tau} =0. 
      $$
      Moreover, $V_p$ is an interaction term which is small with respect to the previous two (in the sense specified in the statement of~\cite[Proposition 3.1]{BiaSpi:ode}).  
\end{lemma}
Note that, in particular, the statement of Lemma~\ref{l:slaving} ensures that any solution lying on $\mathcal M$ becomes exponentially close to an orbit on $\mathcal S_0$ and hence $\mathcal M$ can be regarded as a slaving manifold for $\mathcal S_0$. The proof of Lemma~\ref{l:slaving} is a bit technical, but it can be obtained by relying on a fixed point argument. 

Toward a proof of Theorem~\ref{t:sode:stable}, one first defines $\mathcal S_0$ as follows. Consider the restriction of~\eqref{e:sode3} to the manifold of the slow dynamics: as pointed out before, this system is non singular and $\mathcal E$ is a manifold of equilibria. Hence, one can directly apply the Uniformly Stable Manifold Theorem~\ref{t:us} and define $\mathcal S_0$ as the uniformly stable manifold relative to $\mathcal E$ of system~\eqref{e:sode3} restricted to the manifold of the slow dynamics. By applying Lemma~\ref{l:slaving} and by relying on some more work one can then conclude the proof of Theorem~\ref{t:sode:stable}. 
\end{enumerate}

\bibliography{biblio_trieste}

\end{document}